\documentclass[11pt]{article} 
\usepackage{latexsym} 
\usepackage{amsmath}
\usepackage{amssymb} 
\usepackage{amsfonts}
\usepackage{amscd} 
\usepackage{graphicx}
\usepackage{layout}
\begin{document} 
\newtheorem{Th}{Theorem}[section]
\newtheorem{Cor}{Corollary}[section]
\newtheorem{Prop}{Proposition}[section]
\newtheorem{Lem}{Lemma}[section]
\newtheorem{Def}{Definition}[section]
\newtheorem{Rem}{Remark}[section]
\newtheorem{Ex}{Example}[section]
\newtheorem{stw}{Proposition}[section]


\newcommand{\bet}{\begin{Th}}
\newcommand{\ent}{\stepcounter{Cor}
   \stepcounter{Prop}\stepcounter{Lem}\stepcounter{Def}
   \stepcounter{Rem}\stepcounter{Ex}\end{Th}}


\newcommand{\bec}{\begin{Cor}}
\newcommand{\enc}{\stepcounter{Th}
   \stepcounter{Prop}\stepcounter{Lem}\stepcounter{Def}
   \stepcounter{Rem}\stepcounter{Ex}\end{Cor}}
\newcommand{\bep}{\begin{Prop}}
\newcommand{\enp}{\stepcounter{Th}
   \stepcounter{Cor}\stepcounter{Lem}\stepcounter{Def}
   \stepcounter{Rem}\stepcounter{Ex}\end{Prop}}
\newcommand{\bel}{\begin{Lem}}
\newcommand{\enl}{\stepcounter{Th}
   \stepcounter{Cor}\stepcounter{Prop}\stepcounter{Def}
   \stepcounter{Rem}\stepcounter{Ex}\end{Lem}}
\newcommand{\bef}{\begin{Def}}
\newcommand{\enf}{\stepcounter{Th}
   \stepcounter{Cor}\stepcounter{Prop}\stepcounter{Lem}
   \stepcounter{Rem}\stepcounter{Ex}\end{Def}}
\newcommand{\ber}{\begin{Rem}}
\newcommand{\enr}{
   \stepcounter{Th}\stepcounter{Cor}\stepcounter{Prop}
   \stepcounter{Lem}\stepcounter{Def}\stepcounter{Ex}\end{Rem}}
\newcommand{\bee}{\begin{Ex}}
\newcommand{\ene}{
   \stepcounter{Th}\stepcounter{Cor}\stepcounter{Prop}
   \stepcounter{Lem}\stepcounter{Def}\stepcounter{Rem}\end{Ex}}
\newcommand{\Proof}{\noindent{\it Proof\,}:\ }
\newcommand{\beP}{\Proof}
\newcommand{\enP}{\hfill $\Box$ \par\vspace{5truemm}}

\newcommand{\EE}{\mathbf{E}}
\newcommand{\QQ}{\mathbf{Q}}
\newcommand{\R}{\mathbf{R}}
\newcommand{\C}{\mathbf{C}}
\newcommand{\ZZ}{\mathbf{Z}}
\newcommand{\KK}{\mathbf{K}}
\newcommand{\NN}{\mathbf{N}}
\newcommand{\PP}{\mathbf{P}}
\newcommand{\HH}{\mathbf{H}}
\newcommand{\uuu}{\boldsymbol{u}}
\newcommand{\xxx}{\boldsymbol{x}}
\newcommand{\aaa}{\boldsymbol{a}}
\newcommand{\bbb}{\boldsymbol{b}}
\newcommand{\AAA}{\mathbf{A}}
\newcommand{\BBB}{\mathbf{B}}
\newcommand{\ccc}{\boldsymbol{c}}
\newcommand{\ddd}{\boldsymbol{d}}
\newcommand{\iii}{\boldsymbol{i}}
\newcommand{\jjj}{\boldsymbol{j}}
\newcommand{\kkk}{\boldsymbol{k}}
\newcommand{\rrr}{\boldsymbol{r}}
\newcommand{\FFF}{\boldsymbol{F}}
\newcommand{\yyy}{\boldsymbol{y}}
\newcommand{\ppp}{\boldsymbol{p}}
\newcommand{\qqq}{\boldsymbol{q}}
\newcommand{\nnn}{\boldsymbol{n}}
\newcommand{\vvv}{\boldsymbol{v}}
\newcommand{\eee}{\boldsymbol{e}}
\newcommand{\fff}{\boldsymbol{f}}
\newcommand{\hhh}{\boldsymbol{h}}
\newcommand{\gggg}{\boldsymbol{g}}
\newcommand{\www}{\boldsymbol{w}}
\newcommand{\0}{\boldsymbol{0}}
\newcommand{\lon}{\longrightarrow}
\newcommand{\ga}{\gamma}
\newcommand{\pa}{\partial}
\newcommand{\QED}{\hfill $\Box$}
\newcommand{\id}{{\mbox {\rm id}}}
\newcommand{\Ker}{{\mbox {\rm Ker}}}
\newcommand{\grad}{{\mbox {\rm grad}}}
\newcommand{\ind}{{\mbox {\rm ind}}}
\newcommand{\rot}{{\mbox {\rm rot}}}
\newcommand{\diver}{{\mbox {\rm div}}}
\newcommand{\Gr}{{\mbox {\rm Gr}}}
\newcommand{\Fl}{{\mbox {\rm Fl}}}
\newcommand{\St}{{\mbox {\rm St}}}
\newcommand{\LC}{\mbox{\rm {LC}}}
\newcommand{\LG}{{\mbox {\rm LG}}}
\newcommand{\Diff}{{\mbox {\rm Diff}}}
\newcommand{\Symp}{{\mbox {\rm Symp}}}
\newcommand{\Ct}{{\mbox {\rm Ct}}}
\newcommand{\Uns}{{\mbox {\rm Uns}}}
\newcommand{\rank}{{\mbox {\rm rank}}}
\newcommand{\sign}{{\mbox {\rm sign}}}
\newcommand{\Spin}{{\mbox {\rm Spin}}}
\newcommand{\Sp}{{\mbox {\rm sp}}}
\newcommand{\Int}{{\mbox {\rm Int}}}
\newcommand{\Hom}{{\mbox {\rm Hom}}}
\newcommand{\Tan}{{\mbox {\rm Tan}}}
\newcommand{\codim}{{\mbox {\rm codim}}}
\newcommand{\ord}{{\mbox {\rm ord}}}
\newcommand{\Iso}{{\mbox {\rm Iso}}}
\newcommand{\corank}{{\mbox {\rm corank}}}
\def\mod{{\mbox {\rm mod}}}
\newcommand{\pt}{{\mbox {\rm pt}}}
\newcommand{\qed}{\hfill $\Box$ \par}
\newcommand{\spe}{\vspace{0.4truecm}}
\newcommand{\ad}{{\mbox{\rm ad}}}
\newcommand{\OO}{{\mbox {\rm O}}}
\newcommand{\GL}{{\mbox {\rm GL}}}
\newcommand{\SO}{{\mbox {\rm SO}}}
\newcommand{\SU}{{\mbox {\rm SU}}}
\newcommand{\pr}{{\mbox {\rm pr}}}
\def\mod{{\mbox {\rm \ mod.}}}
\newcommand{\arccosh}{{\mbox{\rm arccosh}}}

\newcommand{\dint}[2]{{\displaystyle\int}_{{\hspace{-1.9truemm}}{#1}}^{#2}}

%
\newenvironment{FRAME}{\begin{trivlist}\item[]
	\hrule
	\hbox to \linewidth\bgroup
		\advance\linewidth by -10pt
		\hsize=\linewidth
		\vrule\hfill
		\vbox\bgroup
			\vskip5pt
			\def\thempfootnote{\arabic{mpfootnote}}
			\begin{minipage}{\linewidth}}{%
			\end{minipage}\vskip5pt
		\egroup\hfill\vrule
	\egroup\hrule
	\end{trivlist}}

\title{
Singularities of Frontals
} 

\author{Goo ISHIKAWA\thanks{This work was supported by JSPS KAKENHI No.15H03615 and No.15K13431.}
}


\date{ }

\maketitle

%
%
%



\section{Introduction}

\

In this survey article we introduce the notion of frontals, which provides a class of 
generalised submanifolds with singularities but with well-defined tangent spaces. 
We present a review of basic theory and 
known studies on frontals in several geometric problems from singularity theory viewpoints. 
In particular, in this paper, we try to give some of 
detailed proofs and related ideas, which were omitted in the original papers, to the basic and important results related to frontals.

We start with one of theoretical motivations for our notion \lq\lq frontal\rq\rq. 
Let $M$ be a $C^\infty$ manifold of dimension $m$, which is regarded as an ambient space. 
Let $f : N \to M$ be an {\it immersion} of an $n$-dimensional $C^\infty$ manifold $N$, which is 
regarded as a parameter space, to $M$. 
Then for each point $t \in N$, we have the $n$-plane $f_*(T_tN)$, the image of the differential map $f_* : T_tN 
\to T_{f(t)}M$ at $t$ in the tangent space $T_{f(t)}M$. Thus we have a field of tangential $n$-planes $\{f_*(T_tN)\}_{t \in N}$ 
along the immersion $f$. 
Moreover if $M$ is endowed with a Riemannian metric, then we have also 
a field of tangential $(m - n)$-planes $f_*(T_tN)^\perp$ along $f$. 
From those vector bundles we can develop differential topology, 
theory of characteristic classes and so on of immersed submanifolds. 
Besides, taking local adapted frames for immersions, we can develop differential geometry of immersed submanifolds in terms 
of frames. Then a natural and challenging problem arises to us on the possibility to find a natural class of
singular mappings enjoying the same properties as immersed submanifolds and 
to develop generalised topological and geometric theories on them. 

In this paper we introduce such a class of generalised submanifolds in terms of Grassmannians: 
Let $\Gr(n, TM)$ denote the Grassmannian of tangential $n$-planes in the tangent bundle $TM$ over 
an $m$-dimensional $C^\infty$ manifold $M$ with the canonical projection $\pi : \Gr(n, TM) \to M$ 
(see \S \ref{Grassmannian bundle and canonical differential system}). 
Let $N$ be a $C^\infty$ manifold of dimension $n$ with $0 \leq n \leq m$ and take a point $a \in N$. 
Then a $C^\infty$ map-germ $f : (N, a) \to M$ is called a {\bf frontal map-germ} or a {\bf frontal} in short 
if there exists a \lq\lq Legendre\rq\rq\, lifting of $f$, that is, there exist an open neighbourhood $U$ of $a$ and 
a $C^\infty$ lift $\widetilde{f} : U \to \Gr(n, TM)$ of $f$, $\pi\circ\widetilde{f} = f\vert_U$, 
such that the image of differential $f_*(T_tN)$ is contained in 
$\widetilde{f}(t)$, for any $t \in U$. Note that $\widetilde{f}(t)$ is an $n$-plane in $T_{f(t)}M$. 
Moreover a $C^\infty$ mapping $f : N \to M$ is called a {\bf frontal mapping} or a {\bf frontal} in short 
if, the germ $f : (N, a) \to M$ at any point $a \in N$ is a frontal. 
See \S \ref{Generalised frontals} for details. 
The formulation using Grassmannians is very natural and satisfactory from the viewpoint of differential systems 
and their geometric solutions as well. See for instance \cite{Yamaguchi08}\cite{IM}\cite{IM2}.

Note that, if $\dim(N) = 1$, then any frontal $f : N \to M$ has 
a global Legendre lift $\widetilde{f} : N \to \Gr(1, TM)$ (Lemma \ref{global-n=1}). 
However, if $\dim(N) =2$, then a frontal $f : N \to M$ not necessarily has a global Legendre lift (Example \ref{global}). 
This fact seems to be found first in the present paper. 
Also note that any mapping $f : N \to M$ is a frontal if $\dim(N) = \dim(M)$ 
(Remark \ref{equi-dimensional}). Any constant mapping $f : N \to M$ is a frontal. 

The notion of \lq\lq frontals\rq\rq\, was introduced already 
in many papers, e.g. \cite{Ishikawa00-2}\cite{ZK}\cite{SUY09-1}\cite{Nuno-Ballesteros}\cite{Brander1}\cite{Brander2}, 
in the case of hypersurfaces as a natural generalisation of wave-fronts. See \S \ref{The case of hypersurfaces}. 

We are going to give a survey on local classification of singularities appearing in frontals in various geometric contexts. 
Basically we mean by the \lq\lq singularities\rq\rq\, of frontals, as usual, the equivalence classes 
of germs of frontals under the following equivalence relation: 

\bef
{\rm 
Two map-germs $f : (N, a) \to (M, f(a))$ and $g : (N', a') \to (M', f'(a'))$ are 
{\bf right-left equivalent} or {\bf ${\mathcal A}$-equivalent} or {\bf diffeomorphic}, 
if there exist diffeomorphism-germs 
$\varphi : (N, a) \to (N', a')$ and $\Phi : (M, f(a)) \to (M', f'(a'))$ 
such that the following diagram commutes: 
$$
\begin{array}{ccc}
(N, a) & \xrightarrow{\  f\ } & (M, f(a)) \\
\varphi\downarrow\hspace{0.2truecm} & & \hspace{0.2truecm}\downarrow \Phi 
\\
(N', a') & \xrightarrow{\ g\ } & (M', f'(a')). 
\end{array}
$$
}
\enf

As the typical singularities of frontals, we introduce cuspidal edges, swallowtails, folded umbrellas, open swallowtails, open folded umbrellas and so on. 

The {\bf cuspidal edge} is defined as the equivalence class of the map-germ $(\R^2, 0) \to (\R^m, 0)$, $m \geq 3$, 
$$
(t, s) \mapsto (t + s, \ t^2 + 2st, \ t^3 + 3st^2, \ 0, \ \dots, \ 0), 
$$
which is diffeomorphic to $(u, w) \mapsto (u, w^2, w^3, 0, \dots, 0)$. 
The cuspidal edge singularities are originally defined only in the three dimensional space. 
Here we are generalising the notion of the cuspidal edge in higher dimensional ambient space. 
It will be often emphasised it by writing \lq\lq embedded" cuspidal edge.

The {\bf folded umbrella} (or the {\bf cuspidal cross cap}) is defined as the equivalence class of the map-germ 
$(\R^2, 0) \to (\R^3, 0)$, 
$$
(t, s) \mapsto (t + s, \ t^2 + 2st, \ t^4 + 4st^3), 
$$
which is diffeomorphic to $(u, t) \mapsto (u, \ t^2 + ut, \ t^4 + {\textstyle \frac{2}{3}}ut^3)$. 

The {\bf open folded umbrella} is defined as the equivalence class of the map-germ $(\R^2, 0) \to (\R^m, 0)$, $m \geq 4$, 
$$
(t, s) \mapsto (t + s, \ t^2 + 2st, \ t^4 + 4st^3, \ t^5 + 5st^4, \ 0, \ \dots, \ 0), 
$$
which is diffeomorphic to $(u, t) \mapsto (u, \ t^2 + ut, \ t^4 + \frac{2}{3}ut^3, \ t^5 + \frac{5}{4}ut^4, \ 0, \ \dots, \ 0)$. 
The open folded umbrella appeared for instance as a frontal-symplectic singularity in the paper \cite{IJ03}. 

The {\bf swallowtail} is defined as the equivalence class of the map-germ $(\R^2, 0) \to (\R^3, 0)$,  
$$
(t, s) \mapsto (t^2 + s, \ t^3 + {\textstyle \frac{3}{2}}st, \ t^4 + 2st^2), 
$$
which is diffeomorphic to $(u, t) \mapsto (u, \ t^3 + ut, \ t^4 + \frac{2}{3}ut^2)$. 

The {\bf open swallowtail} is defined as the equivalence class of the map-germ $(\R^2, 0) \to (\R^m, 0)$, $m \geq 4$, 
$$
(t, s) \mapsto (t^2 + s, \ t^3 + {\textstyle \frac{3}{2}}st, \ t^4 + 2st^2, \ t^5 + {\textstyle \frac{5}{2}}st^3, \ 0, \ \dots, \ 0), 
$$
which is diffeomorphic to $(u, t) \mapsto (u, \ t^3 + ut, t^4 + \frac{2}{3}ut^2, \ t^5 + \frac{5}{9}ut^3, \ 0, \ \dots, \ 0)$. 
The open swallowtail singularity was introduced by Arnol'd (see \cite{Arnold81}) as a singularity of Lagrangian varieties in symplectic geometry. Here we abstract its diffeomorphism class as the singularity of parametrised surfaces (see \cite{Givental86}\cite{Ishikawa12-2}). 


\

In Part I, we provide basic studies for an intrinsic understanding of frontals 
as parametrised singular submanifolds with well-defined tangent spaces. 

We give the exact definition of frontals in \S \ref{The case of hypersurfaces} in the case of hypersurfaces and, after the description of 
Grassmannian bundles and canonical (or generalised contact) distributions 
in \S \ref{Grassmannian bundle and canonical differential system}, we give the general definition in \S \ref{Generalised frontals}. 
In \S \ref{Density function}, we have introduced the density function as a main notion for the theory of frontals. 

A frontal $f : N \to M$ is called a {\bf proper frontal} in the present paper if the singular (non-immersive) locus $S(f)$ is 
nowhere dense in $N$ (\S \ref{Proper frontals}). 
In \cite{Ishikawa13}\cite{Ishikawa14}\cite{Ishikawa15}\cite{Ishikawa15-2}, \lq\lq frontal\rq\rq\, maps were defined as proper frontals, namely, the density of regular locus was assumed. 
Note that proper frontals are not generic in the space of all frontals 
for $C^\infty$-topology in general (Remark \ref{non-approximation-}). 
In \S \ref{tangent-and-normal}, we introduce the tangent bundles to frontals. 

Viewed from our generalisation, the notion of frontals turns to be closely related to the notion of 
{\bf openings}. Though the notion of openings of mappings 
seems to be noticed naively in many previous contributions, 
it is introduced in the author's recent papers \cite{Ishikawa13}\cite{Ishikawa14}. 
An opening separates the self-intersections of the original map-germ, preserving its singularities.  
For example, the swallowtail is an opening of the {\bf Whitney's cusp} map-germ $(\R^2, 0) \to (\R^2, 0)$ 
defined by $(t, s) \mapsto (t^2 + s, \ t^3 + {\textstyle \frac{3}{2}}st)$ which is diffeomorphic 
to $(u, t) \mapsto (u, \ t^3 + ut)$ and 
the open swallowtail is a versal opening of them. 
Openings of map-germs appear as typical singularities in several problems of geometry and its applications. Note that the process of unfoldings of map-germs $(\R^n, 0) \to (\R^m, 0)$ 
preserves the \lq\lq relative dimension\rq\rq\, $m - n$. 
On the other hand, the process of openings preserves $n$ but changes $m$, and it 
gives bridges between map-germs of different relative dimensions. 
We recall also the related notions, \lq\lq Jacobi modules\rq\rq\, ${\mathcal J}_f$ 
and \lq\lq ramification modules\rq\rq\,  ${\mathcal R}_f$. They play important role 
to analysis and classification of singularities of mappings $f$,  
in particular, the study on symplectic singularities, 
contact singularities and singularities of tangent surfaces 
(\cite{Ishikawa83}\cite{Ishikawa92}\cite{Ishikawa92-2}\cite{Ishikawa95}\cite{Ishikawa96}\cite{Ishikawa00}\cite{Ishikawa04}\cite{Ishikawa05}\cite{IJ08}). Moreover those notions seem to be related 
to recognition problem of singularities (see Definition \ref{J-equivalence}). 
Note that we used the notation, for the ramification module of $f$, \lq${\mathcal D}_f$\rq \ instead of ${\mathcal R}_f$ in \cite{Ishikawa83}, 
relating Mather's ${\mathcal C}$-equivalence, and 
we denoted it by \lq$H_f$\rq \ in \cite{Ishikawa92}\cite{Ishikawa92-2}\cite{Ishikawa94}\cite{Ishikawa95}, 
because it can be regarded as a cohomological invariant. 
Note that the notion of openings, Jacobi modules and ramification modules for {\it multi-germs} is naturally introduced in the paper \cite{Ishikawa13}. We give a review on the theory of opening related to frontals in \S \ref{Openings and frontals} and \S \ref{Versal openings}. 
Moreover in \S \ref{Subfrontals and superfrontals} we give ideas of \lq\lq subfrontals\rq\rq\, and \lq\lq superfrontals\rq\rq, related to frontals.

In \cite{SUY12}, it is introduced the related notion of \lq\lq coherent tangent bundles\rq\rq 
as generalised Riemannian manifolds. 
Moreover Saji, Umehara, Yamada are developing the intrinsic studies of frontals in terms of singular metrics introduced by Kossowski \cite{Kossowski04}. 
We intend to give abstract differential-topological features of frontals, which is invariant under diffeomorphisms, 
by proving another way to study intrinsically frontals in terms of the theory of $C^\infty$-rings (\S \ref{Algebraic openings}). 

In part II, we give a survey of several results on frontals as an application of the basic theory 
presented in part I. 

In \S \ref{Frontal curves}, we treat frontal curves and give basic results on them. 
Let $\gamma : N \to \R^2$ a planar frontal curve with $\dim(N) = 1$. 
By Lemma \ref{global-n=1}, 
there exists a global Legendre lifting $\widetilde{\gamma} : I \to P(T^*\R^2)$. 
Thus it is possible to perform differential geometry of planar frontals, as a generalised differential geometry of planar immersions, in terms of Legendre curves covering frontals. 
In fact geometric studies on planar frontals, evolutes and involutes, are given in a series of papers \cite{FT1}\cite{FT2}\cite{FT3}\cite{FT4}. 
As a related topics to planar frontals, 
we gave a review the \lq\lq Goursat Monster tower\rq\rq\, 
found by Zhitomirskii, Montgomery, Mormul and others (cf. \cite{MZ1}\cite{MZ2}\cite{CM}) 
and the \lq\lq Legendre-Goursat duality\rq\rq\, related to it in \cite{Ishikawa15}. 

Singularities and bifurcations of wavefronts based on Legendre singularity theory are established by 
Arnold-Zakalyukin's theory (\cite{Arnold76}\cite{Arnold90}\cite{Zakalyukin1}\cite{Zakalyukin2}). 
The application of singularity theory to differential geometry has been developed by many authors 
(see for instance \cite{Porteous}\cite{Porteous2}\cite{BG}\cite{IFRT}). 
The geometric study of submanifolds in hyperbolic space $H^{n+1}$ based on singularity theory was initiated by Izumiya et al.(\cite{IPS}\cite{IPS2}\cite{IPT}). The Legendre duality developed in 
\cite{Izumiya}\cite{Chen-Izumiya} enables us to unify the theory of framed curves in any space form as describes in \cite{Ishikawa12}. 
We recall Legendre duality (see \cite{Bruce}\cite{Morimoto}\cite{IMo}\cite{IM}\cite{Chen-Izumiya}) in the framework of moving frames and flags and discuss its generalisation and relation with the theory of frontals in \S \ref{Frames and flags}, \S \ref{Legendre duality.} and 
\S \ref{Grassmannian frontals}. 

Let $\gamma : I \to \R^3$ be a space frontal curve. 
Then the tangent surface (tangent developable) $\Tan(\gamma)$ is defined as the surface 
ruled by tangent lines to $\gamma$. Then the tangent surface 
has zero Gaussian curvature, therefore it is flat with respect to Euclidean 
metric of $\R^3$ at least off the singular locus. Thus the tangent surfaces serve 
main parts of \lq\lq flat frontals\rq\rq (\S \ref{Tangent varieties}). 
Flat fronts or flat frontals are studied also in \cite{MU}\cite{Naokawa}. 

The notion of tangent surfaces ruled by \lq\lq tangent lines\rq\rq\ to directed curves 
is naturally generalised in various ways: For a curve in a projective space, 
we regard tangent projective lines as \lq\lq tangent lines\rq\rq. 
The classification is generalised to $A_n$-geometry (\S \ref{Grassmannian geometry}). 
For a curve in a Riemannian manifold, 
we regard tangent geodesics as \lq\lq tangent lines. 
In fact, tangent surfaces are defined for proper frontal curves (directed curves) 
in a manifold with an affine connection (\S \ref{Affine connection and tangent surface}). 
After discussing useful criteria of singularities in \S \ref{Characterisation of frontal singularities}, we define 
null tangent surface 
to a null curve of a semi (pseudo)-Riemannian manifold, 
regarding null geodesics as \lq\lq tangent lines\rq\rq (see \S \ref{Null frontals}). 
In particular we pick up several results related to $D_n$-geometry (\cite{IMT4}). 
For a horizontal curve of a sub-Riemannian manifold, 
we regard \lq\lq tangent lines\rq\rq\, by abnormal geodesics (see \S \ref{Abnormal frontals}). 
In particular the classification result of singularities of tangent surfaces to generic 
integral curves to Cartan distribution with $G_2$-symmetry is introduced. 

Speaking of $G_2$, we note that the work on frontals may be related to the rolling ball problem \cite{Agrachev}\cite{BM}\cite{BH}\cite{Montgomery}. 
We will treat \lq\lq rolling frontals\rq\rq\, as a generalisation of rolling bodies \cite{AN} in a forthcoming paper. 

In the last section (\S \ref{Appendix: Malgrange preparation theorem on differentiable algebras}), 
as an appendix, we show the Malgrange's preparation theorem on differentiable algebras (\cite{Malgrange}) 
from the ordinary Malgrange-Mather's preparation theorem (see for example \cite{Brocker}), relating to the theory of 
$C^\infty$-rings which we have utilised in this paper. 

The author hopes very much that this survey paper helps to raise wider reader's interest to the mathematics on frontals. 

\

In this paper a manifold or a mapping is supposed to be of class $C^\infty$ 
unless otherwise stated. The symbol $\subseteq$ of inclusion is often used, which 
has the same meaning as $\subset$
 just to stress that the equality may occur. 

\

\

\begin{center}
{\bf {\LARGE Part I. Basic Theory}}
\end{center}

\section{The case of hypersurfaces}
\label{The case of hypersurfaces}

Let $M$ be a manifold of dimension $m$. 
Let $P(T^*M)$ denote the projective cotangent bundle of $M$, 
which consists of non-zero cotangent vectors somewhere on $M$ 
considered up to a non-zero scalar multiplication. 
Note that $P(T^*M)$ is naturally identified with the Grassmannian bundle $\Gr(m-1, TM)$ (see \S \ref{Grassmannian bundle and canonical differential system}) by sending 
each class $(x, [\alpha]) \in P(T^*M)$ of a non-zero covector $\alpha \in T^*M$ to 
its kernel $\Ker(\alpha) \in \Gr(m-1, T_xM)$. 
Note that $\alpha \in T_x^*M$, that $\alpha : T_xM \to \R$ is a non-zero linear map, 
and that $\Ker(\alpha) \subset T_xM$ is an $(m-1)$-plane. 
Then the $(2m - 1)$-dimensional manifold 
$P(T^*M)$ has a canonical contact structure $D \subset TP(T^*M)$. 
In fact it is defined by $D = \bigcup_{(x, [\alpha])} D_{(x, [\alpha])}$, 
and $D_{(x, [\alpha])} = \pi_*^{-1}(\Ker(\alpha))$, where 
$\pi : P(T^*M) \to M$ is the canonical projection. 

We recall the coordinate description of the contact structure, which will be needed for the detailed 
computation on singularities. 

Let $(x^1, x^2, \dots, x^m)$ be a local coordinate system on an open subset $U$ of $M$. 
Let 
$$(x^1, x^2, \dots, x^m, p_1, p_2, \dots, p_m)
$$ 
be the associated system of coordinates on $T^*U$ 
such that any element $\alpha \in T^*U$ is expressed as
$$
\alpha = p_1dx^1 + p_2dx^2 + \cdots + p_mdx^m, 
$$
by its coordinates. Set $V_i = \{ p_i \not= 0\} \subset T^*U, 1 \leq i \leq m$. 
Then we have a local system of coordinates of $P(T^*M)$ associated to $V_i$, 
$$
x^1, x^2, \dots, x^m, - p_1/p_i, \dots, - p_{i-1}/p_i, - p_{i+1}/p_i, \dots, - p_m/p_i. 
$$
To avoid non-essential complexity, we will discuss just for $i = m$ in what follows. 
Then set $a_i = - p_i/p_m, 1 \leq i \leq m-1$. Then 
$$
x^1, x^2, \dots, x^m, a_1, a_2, \dots, a_{m-1}
$$
give a local system of coordinates of $P(T^*M)$ 
and the contact structure $D \subset TP(T^*M)$ is given locally by 
$$
dx^m - (a_1dx^1 + a_2dx^2 + \cdots + a_{m-1}dx^{m-1}) = 0. 
$$

Let $N$ be a submanifold of dimension $n$ with $n < m$. 
Then the submanifold $N$ induces the {\bf projective conormal bundle} 
$$
\widetilde{N} = P(T^*_NM) = \{ (x, [\alpha]) \in P(T^*M) 
\mid \alpha\vert_{T_xN} = 0 \}, 
$$
which satisfies that $T\widetilde{N} \subset D$ and $\dim(\widetilde{N}) = m-1$, 
in other words, 
a Legendre submanifold in the contact manifold $P(T^*M)$.  

In particular, suppose $n = m - 1$, that is, $N$ is a hypersurface of $M$. 
Then $\pi\vert_{\widetilde{N}} : \widetilde{N} \to N$ is a diffeomorphism. 
Its inverse $N \to \widetilde{N}$ is given by $x \mapsto (x, T_xN)$. 

Let $f : N \to M$ be an immersion of an $(m - 1)$-dimensional manifold $N$ to 
an $m$-dimensional manifold $M$. 
Then we have an immersion $\widetilde{f} : N \to P(T^*M)$ 
defined by $\widetilde{f}(t) = (f(t), f_*(T_tN))$. 
Then $\widetilde{f}$ is a lift of $f$ and $\widetilde{f}$ is $D$-integral, i.e. 
$\widetilde{f}_*(T_tN) \subset D_{f(t)}$ for any $t \in N$. 
In other words, $\widetilde{f}$ is a Legendre immersion.

\ber
{\rm
Set $\widetilde{f}(t) = (f(t), [\alpha(t)]) \in P(T^*_{f(t)}M)$. Then 
the condition $\widetilde{f}_*(T_tN) \subset D_{f(t)}$ is equivalent to that 
$\alpha(t)\vert_{f_*(T_tN)} = 0$. 
}
\enr

\bef
{\rm
Let $N$ be a manifold of dimension $m - 1$ in a manifold $M$ of dimension $m$. 
A map-germ $f : (N, a) \to M$ is called a  {\bf wave-front} or a {\bf front} in short 
if there exists a germ of Legendre immersion $\widetilde{f} : (N, a) \to 
P(T^*M)$ with $\pi\circ\widetilde{f} = f$. 

A mapping $f : N \to M$ is called a {\bf wave-front} or a {\bf front} in short 
if, for any point $a \in N$, the germ of $f$ at $a$ is a front. 
}
\enf

A map-germ $f : (N, a) \to M$ is a front 
if and only if there exists a representative of $f$, which is a front. 

\ber
{\rm
In the original and naive context, the image $f(N)$ was called a wave-front rather than the parametrisation $f$ itself. 
}
\enr

\bef
\label{frontal-hypersurface2}
{\rm
Let $N$ be a manifold of dimension $m - 1$ in a manifold $M$ of dimension $m$. 
Let $a \in N$. 
A map-germ $f : (N, a) \to M$ is called a  {\bf frontal map-germ} or a {\bf frontal} in short 
if there exist a germ of Legendre lifting 
$\widetilde{f} : (N, a) \to P(T^*M)$ of $f$, that is, 
there exist an open neighbourhood 
$U$ of $a$, a representative $f : U \to M$ of $f$ and 
a Legendre lifting $\widetilde{f} : U \to P(T^*M)$ of $f\vert_U$, i.e. 
$\widetilde{f}_*(T_tN) \subset D_{f(t)}$ for any $t \in U$ and 
$\pi\circ \widetilde{f} = f\vert_U$. 
Here we do not assume that $\widetilde{f}$ is an immersion. 

A mapping $f : N \to M$  is called a {\bf frontal mapping} or a {\bf frontal} in short 
if, for any $a \in N$, 
the germ of $f$ at $a$ is a frontal. 
}
\enf

A map-germ $f : (N, a) \to M$ is a frontal if and only if 
there exists a representative of $f$, which is a frontal. 

\

In Definition \ref{frontal-hypersurface2} we have defined the notion of frontals by the {\it local} 
existence of its Legendre liftings. A frontal $f : N \to M$ 
not necessarily has its {\it global} Legendre lifting $\widetilde{f} : N \to P(T^*M)$. 

\bee
\label{global}
{\rm
Define a $C^\infty$ function $\varphi : \R \to \R$ by $\varphi(t) = e^{-1/t^2} (t > 0), \varphi(t) = 0 (t \leq 0)$. 
Then define $h : \R^2 \to \R^3$ by $h(t_1, t_2) = (t_1, t_2^2, t_2^3 + \varphi(t_1)t_2)$, which we will call 
a {\bf half cuspidal edge}. 

\


\

The mapping $h$ is not frontal. In fact the local existence of Legendre lift for $h$ does not 
hold at the origin $(t_1, t_2) = (0, 0)$. Moreover $h$ is a front on $\R^2 \setminus \{(0, 0)\}$ 
with cuspidal edge along $\{ (t_1, 0) \mid t_1 < 0\}$ 
and the Legendre lifting $\widetilde{h} : \R^2 \setminus \{(0, 0)\} \to P(T^*\R^3) \cong \R^3\times \R P^2$ 
is not homotopically trivial. In fact $\widetilde{h}$ restricted to a loop around the origin of $\R^2$ 
generates the fundamental group $\pi_1(P(T^*\R^3)) \cong \pi_1(\R P^2) \cong \ZZ/2\ZZ$. 

Define $k : \R^2 \to \R^2$ by $k(t_1, t_2) = \varphi(t_1^2 + t_2^2 - 1)(t_1, t_2)$. 
Then $k$ is a $C^\infty$ mapping which collapses the unit disc to the origin and 
maps the outside of the unit disc to $\R^2 \setminus \{(0, 0)\}$ diffeomorphically. 
Set $f = h\circ k : \R^2 \to \R^3$. Then 
\\
(1) $f$ is a frontal.
\\
(2) There does not exist a global Legendre lifting $\widetilde{f} : \R^2 \to P(T^*\R^3)$ of $f$. 
\\
To see (1), let $t = (t_1, t_2)$ satisfy $t_1^2 + t_2^2 < 1$. Then the germ of $f$ at $t$ is constant and therefore 
it is a frontal. Let $t$ satisfy $t_1^2 + t_2^2 > 1$. Then the germ of $f$ at $t$ is right equivalent to 
$h$ at $k(t) \in \R^2 \setminus \{(0, 0)\}$, which is a frontal. Let $t$ satisfy $t_1^2 + t_2^2 = 1$. 
Then any local extension of $\widetilde{h}\circ k$ to $(\R^2, t)$ turns to be a Legendre lift of the germ of $f$ at $t$. 
Therefore $f$ is a frontal. Thus we have (1). To see (2), it is sufficient to observe that $\widetilde{h}\circ k : \R^2 \setminus D^2 
\to P(T^*\R^3)$, which is the unique Legendre lift of $f$ restricted to $\R^2 \setminus D^2$,  
is never extended continuously to $\R^2$. 
}
\ene

\section{Grassmannian bundle and generalised contact distribution}
\label{Grassmannian bundle and canonical differential system}

Let $N$ be a manifold of dimension $n$ and $M$ a manifold of dimension $m$ with $n \leq m$. 
Note that in the previous section we assumed $m = n+1$. However in the next section 
we treat the general case under the weaker condition $n \leq m$. 

To treat the general cases, we recall the Grassmannian bundle associated to the tangent bundle 
$TM$ of $M$. 
For each $x \in M$, $\Gr(n, T_xM)$ denotes the Grassmannian manifold consisting of 
$n$-dimensional subspaces 
of $T_xM$. 
Then let $\Gr(n, TM) = \bigcup_{x \in M} \Gr(n, T_xM)$. 
Note that $\Gr(n, TM)$ is a bundle over $M$ with fibres $\Gr(n, T_xM)$ and 
that $\dim(\Gr(n, T_xM)) = n(m - n)$. 
Note also that $\Gr(n, T_xM)$ is identified with $\Gr(m-n, T^*_xM)$ 
and therefore that, when $m = n+1$, $\Gr(n, TM)$ is identified with $P(T^*M)$. 
Let $\pi : \Gr(n, TM) \to M$ be the canonical projection, $\pi(x, V) = x$ for any $(x, V) \in \Gr(n, TM)$ with 
$V \in \Gr(n, T_xM), x \in M$. 
If $n = m$, then $\pi : \Gr(m, TM) \to M$ is a diffeomorphism. 

\bel
\label{canonical--canonical}
Let $\Phi : M \to M'$ be a diffeomorphism. Let $n$ be an integer with 
$0 \leq n \leq m = \dim(M)$. 
Let $\Phi_\sharp : \Gr(n, TM) \to \Gr(n, TM')$ 
denote the diffeomorphism induced by the differential map $\Phi_*$ which is regarded as the 
bundle isomorphism $\Phi_* : TM \to TM'$ covering $\Phi$. 
Then we have 
$\pi\circ \Phi_{\sharp} = \Phi\circ\pi : \Gr(n, TM) \to M'$. 
Here $\pi$ means the canonical projection $\Gr(n, TM') \to M'$ as well as $\Gr(n, TM) \to M$. 
\enl

\Proof
Let $(x, V) \in \Gr(n, TM)$. Then $\Phi_{\sharp}(x, V) = (\Phi(x), \Phi_*(V))$. Therefore 
$
(\pi\circ \Phi_{\sharp})(x, V) = \pi(\Phi(x), \Phi_*(V)) = \Phi(x) = (\Phi\circ\pi)(x, V). 
$
\QED

\

We recall the coordinate description of Grassmannians. 
Let $(x_0, V_0) \in  \Gr(n, TM)$. 
Here $x_0 \in M$ and $V_0 \in \Gr(n, T_{x_0}M)$ so that $V_0 \subset T_{x_0}M$ is 
a fixed $n$-plane. The suffix $0$ is used to 
indicate that $(x_0, V_0)$ becomes the center of the local coordinate system 
we are going to provide. 
Let us take a local coordinate system $(x^1, \dots, x^n, x^{n+1}, \dots, x^m)$ on a coordinate 
neighbourhood $U \subset M$ 
centred at $x_0$ such that $\pa/\pa x^1, \dots, \pa/\pa x^n$ generate $V_0$ at $x_0$. 
Let $\pi' : U \to \R^n$ denote the coordinate projection defined by $(x^1, \dots, x^n, \dots, x^m) 
\mapsto (x^1, \dots, x^n)$. 
Let $\Omega \subset \pi^{-1}(U)$ be the set of $(x, V)$ with $x \in U, V \in \Gr(n, T_xM)$ 
such that $V$ is mapped isomorphically by $\pi'_* : TU \to T\R^n$ to $\pi'_*(V)$. 
Then, for any $(x, V)$, there exist unique real numbers $a_j^k, (1 \leq j \leq n, n+1 \leq k \leq m)$ 
such that the $n$-plane $V$ has the basis of the form (*)
$$
\left\{
\begin{array}{cccccccc}
h_1 & = & \dfrac{\pa}{\pa x^1}(x) &  & &  & + & a_1^{n+1}\dfrac{\pa}{\pa x^{n+1}}(x) + \cdots + 
a_1^{m}\dfrac{\pa}{\pa x^{m}}(x), 
\vspace{0.3truecm}
\\
h_2 & = & & \dfrac{\pa}{\pa x^2}(x) & & & +  & a_2^{n+1}\dfrac{\pa}{\pa x^{n+1}}(x) + \cdots + 
a_2^{m}\dfrac{\pa}{\pa x^{m}}(x), 
\vspace{0.3truecm}
\\
\vdots &  &  &   & \ddots  & & & 
\vspace{0.3truecm}
\\
h_n & = & & & & \dfrac{\pa}{\pa x^n}(x) & +  & a_2^{n+1}\dfrac{\pa}{\pa x^{n+1}}(x) + \cdots + 
a_2^{m}\dfrac{\pa}{\pa x^{m}}(x). 
\end{array}
\right.
$$
Thus we have a system of coordinates $(x^1, \dots, x^m, a_j^k, (1 \leq j \leq n, n+1 \leq k \leq m))$ 
on $\Omega$ of $\Gr(n, TM)$ centred at $(x_0, V_0)$. 

We call the coordinate systems constructed as above {\bf Grassmannian coordinates}. 

\

The {\bf canonical distribution} $D \subset T(\Gr(n, TM))$ on the Grassmann bundle $\Gr(n, TM)$ 
is defined by $D = \bigcup_{(x, V)} D_{(x, V)}$ where $(x, V)$ runs over $\Gr(n, TM)$, 
$V$ being an $n$-plane of $T_xM$, $x \in M$, 
and, for $v \in T_{(x, V)}(\Gr(n, TM))$, 
$$
v \in D_{(x, V)} \Longleftrightarrow \pi_*(v) \in V (\subset T_xM). 
$$
We call the canonical distribution $D$ on $\Gr(n, TM)$ also the {\bf canonical differential system} and 
also the {\bf contact distribution}, in a generalised and wider sense. If $n = m-1$, then $D$ is 
the contact distribution in the strict sense. 
Note that, if $n = m$, then $D = T(\Gr(n, TM)) \cong TM$. 

\bef
{\rm
A mapping $F : N \to \Gr(n, TM)$ is called an {\bf integral mapping} of the contact distribution $D \subset T\Gr(n, TM)$ 
or a {\bf $D$-integral mapping} if $F_*(TN) \subset D$. 
If $\dim(N) = n$, then we call an integral mapping $f : N^n \to \Gr(n, TM)$ of the contact structure $D \subset T\Gr(n, TM)$ 
a {\bf Legendre mapping} in a generalised and wider sense. 
}
\enf

\bel
\label{D-integral}
A mapping $F : N^n \to \Gr(n, TM)$ is a Legendre mapping if and only if, 
$(\pi\circ F)_*(T_tN) \subseteq F(t), (t \in N)$. 
If $F$ is Legendre and $\pi\circ F$ is an immersion at $t \in N$, 
then $F(t) = (\pi\circ F)_*(T_tN)$. 
\enl

\Proof
By definition, $F$ is Legendre if and only if, for any $t \in N$, $F_*(T_tN) \subset D_{F(t)}$. 
Since $D_{F(t)} = \pi_*^{-1}(F(t))$, regarding $F(t)$ as an $n$-plane in $T_{(\pi\circ F)(t)}M$, 
the condition is equivalent to that 
$\pi_*(F_*(T_tN)) \subseteq F(t)$, that is, $(\pi\circ F)_*(T_tN) \subseteq F(t)$, for any $t \in N$. 
Moreover if $\pi\circ F$ is an immersion at $t \in N$, then $\dim((\pi\circ F)_*(T_tN)) = n$. 
Therefore we have $(\pi\circ F)_*(T_tN) = F(t)$. 
\QED

\

The following result shows one of fundamental properties of the canonical differential systems 
(the generalised contact distributions). 

\bep
\label{canonical-is-canonical}
Let $\Phi : M \to M'$ be a diffeomorphism. Let $0 \leq n \leq m = \dim(M)$. 
Let $D$ denote the contact distribution of $\Gr(n, TM')$ as well as that of $\Gr(n, TM)$. 
Then, for any $(x, V) \in \Gr(n, TM)$, we have
$$
({\Phi_\sharp})_*(D_{(x, V)}) = D_{(\Phi(x), \Phi_*(V))}. 
$$
In particular we have $({\Phi_\sharp})_*(D) = D \subset T(\Gr(n, TM'))$ (see Lemma \ref{canonical--canonical}). 
\enp

\Proof
Let $v \in D_{(x, V)}$. Then $\pi_*v \in V$. Then we have, by Lemma \ref{canonical--canonical}, 
$$
\pi_*((\Phi_{\sharp})_*(v)) = (\pi\circ\Phi_{\sharp})_*(v) = (\Phi\circ\pi)_*(v) 
= \Phi_*(\pi_*v) \in \Phi_*(V). 
$$
Therefore we have 
$({\Phi_\sharp})_*(D_{(x, V)}) \subseteq D_{(\Phi(x), \Phi_*(V))}$. The converse inclusion 
is obtained by considering $\Phi^{-1}$, or, by counting the dimension of the vector spaces. 
\QED

\

We conclude this section by the coordinate description of the contact distribution: 
Take the Grassmannian coordinates $(x_1, \dots, x_m, a_j^k, (1 \leq j \leq n, n+1 \leq k \leq m))$ of $\Gr(n, TM)$ 
on an open set $\Omega \subset \Gr(n, TM)$. 
Set 
$$
\theta^k := dx^k - \sum_{j=1}^n a_j^k dx^j, \ (n+1 \leq k \leq m). 
$$

\bel
\label{canonical-coordinate}
Let $0 \leq n \leq \dim(M)$. 
The local description of the contact distribution $D$ of $\Gr(n, TM)$ is given by 
$$
D\vert_{T\Omega} = \{ v \in T\Omega \mid \theta^{n+1}(v) = 0, \dots, \theta^m(v) = 0 \}. 
$$
\enl

\Proof
Let $(x, V) \in \Omega$ and 
$v = \sum_{i=1}^m b^i\pa/\pa x^i + \sum_{j,k} c_j^k \pa/\pa a_j^k \in T_{(x, V)}\Omega$. 
Then $v \in D_{(x, V)}$ if and only if 
$\pi_*(v) \in V$. Now $V = \langle h_1, \dots, h_n\rangle_\R$ in terms of the above basis (*). Then the condition is equivalent to that 
$\sum_{i=1}^m b^i \pa/\pa x^i = \sum_{j=1}^n \lambda^j h_j$ for some $\lambda^1, \dots, \lambda^n \in \R$, 
which is equivalent to that $b^j = \lambda^j, 1 \leq j \leq n$ and $b^k = \sum_{j=1}^n b^ja_j^k, n+1 \leq k \leq m$, and 
thus equivalent to that $\theta^{k}(v) = 0, n+1 \leq k \leq m$. 
\QED

\section{Generalised frontals}
\label{Generalised frontals}

We give the exact definition of our main notion in this paper: 

\bef
\label{generalised-frontals}
{\rm 
Let $N$ be an $n$-dimensional manifold and $M$ an $m$-dimensional manifold 
with $n \leq m$. 
A map-germ $f : (N, a) \to M$ is called a  {\bf frontal map-germ} or a {\bf frontal} in a generalised sense, if there exists a germ of 
Legendre lift $\widetilde{f} : (N, a) \to \Gr(n, TM)$ of $f$, that is, 
if there exists an open neighbourhood $U$ of $a$ and 
a $D$-integral lift $\widetilde{f} : U \to \Gr(n, TM)$ of $f$ for 
the canonical distribution $D \subset T\Gr(n, TM)$ and for 
the canonical projection $\pi : \Gr(n, TM) \to M$, which 
satisfies that $f_*(T_tN) \subseteq \widetilde{f}(t)$ for any $t \in U$ and 
$\pi\circ \widetilde{f} = f$. 

We call a mapping $f : N^n \to M^m$ a {\bf frontal mapping} or a {\bf frontal} in a generalised sense, 
if, for any point $a \in N$, the germ of $f$ at $a$ is a frontal. 
}
\enf

\ber
\label{equi-dimensional}
{\rm
Note that, in the equi-dimensional case $n = m$, any mapping $f : N \to M$ is a frontal. 
In fact the mapping $\widetilde{f} : N \to \Gr(m, TM)$ defined by $\widetilde{f}(t) := T_{f(t)}M$ 
is a Legendre lift of $f$. 
}
\enr

\bep
\label{induced-D-integral}
Let $f : (N, a) \to (M, f(a))$ and $g : (N', a') \to (M', f(a'))$ be map-germs. 
If $f$ is a frontal and $g$ is right-left equivalent to $f$, then $g$ is a frontal. 
\enp

\Proof
Suppose $g\circ \varphi = \Phi\circ f$ for some diffeomorphism-germs $\varphi : (N, a) \to 
(N', a')$ and $\Phi : (M, f(a)) \to (M', f(a'))$. 
Let $\widetilde{f} : (N, a) \to \Gr(n, TM)$ be a Legendre lift of $f$. 
Set $\widetilde{g} := \Phi_\sharp \circ \widetilde{f} \circ \varphi^{-1} : (N', a') 
\to \Gr(n, TM')$. 
For $t' \in (N, a')$, we have, by Proposition \ref{canonical-is-canonical}, 
$$
{\widetilde{g}}_*(T_{t'}N') = (\Phi_\sharp)_*(\widetilde{f}_*(\varphi^{-1}_*(T_{t'}N)))
= (\Phi_\sharp)_*(\widetilde{f}_*(T_{\varphi^{-1}(t)}N)) 
\subset (\Phi_\sharp)_*D = D. 
$$
Therefore 
$\widetilde{g}$ is Legendre. Moreover, by Lemma \ref{canonical--canonical}, we have 
$$
\pi\circ\widetilde{g} = \pi\circ\Phi_\sharp \circ \widetilde{f} \circ \varphi^{-1} 
= \Phi\circ \pi\circ \widetilde{f} \circ \varphi^{-1} = 
\Phi\circ f \circ \varphi^{-1} = g. 
$$
Therefore $\widetilde{g}$ is a Legendre lifting of $g$, and hence $g$ is a frontal. 
\QED

\

\bef
{\rm 
A map-germ $f : (N, a) \to M$ is called a {\bf front} in the generalised sense 
if there exists a Legendre lift $\widetilde{f} : (N, a) \to \Gr(n, TM)$ of $f$ 
such that $\widetilde{f}$ is an immersion-germ. 
A mapping $f : N^n \to M^m$ is called a {\bf front} in the generalised sense if, 
for any $a \in N$, the germ of $f$ at $a$ is a front. 
}
\enf

A map-germ $f : (N, a) \to M$ is a front in the generalised sense if and only if 
there exists a representative of $f$ which is a front. 
The condition that $ f : N \to M$ is a front in the generalised sense 
is equivalent to the local existence, at each point of $N$, of an immersive lift 
$\widetilde{f} : U \to \Gr(n, TM)$ of $f$ satisfying $f_*(T_tN) \subset \widetilde{f}(t), (t \in U)$.

\section{Density function}
\label{Density function}

The notion of density functions is a key to understand the geometry of frontals, which was introduced in 
\cite{KRSUY}\cite{FSUY}\cite{SUY09-1} first. We introduce its generalisation (see also \cite{IY}\cite{IY1}): 

\bep
\label{density-function}
Let $f : (N, a) \to M$ be a map-germ with $\dim(N) = n \leq m = \dim(M)$. 
Then the following conditions are equivalent:
\\
{\rm (1)} $f$ is a frontal map-germ. 
\\
{\rm (2)} 
There exists a frame $h_1, h_2, \dots, h_n : (N, a) \to TM$ 
along $f$ and a function-germ $\sigma : (N, a) \to \R$ such that 
$$
(\dfrac{\pa f}{\pa t_1} \wedge \dfrac{\pa f}{\pa t_2} \wedge \cdots \wedge \dfrac{\pa f}{\pa t_n})(t) 
= \sigma(t)(h_1 \wedge h_2 \wedge \cdots \wedge h_n)(t), 
$$
as germs of $n$-vector fields $(N, a) \to \wedge^n TM$ over $f$. 
Here $t_1, t_2, \dots, t_n$ are coordinates on $(N, a)$. 
\enp

The function $\sigma : (N, a) \to \R$ in Proposition \ref{density-function} is called 
a {\bf signed area density function} or briefly an {\bf $s$-function} of the frontal $f$ 
associated with the frame. 
Note that the function $\sigma$ is essentially the same thing with the function $\lambda$ introduced in 
\cite{KRSUY}\cite{FSUY} in the case $\dim(M) = 3$. 

Two function-germs $\sigma, \widetilde{\sigma} : (N, a) \to \R$ 
are called {\bf ${\mathcal K}$-equivalent} if 
there exists a diffeomorphism-germ $T : (N, a) \to (N, a)$ 
and a non-vanishing function-germ $c : (N, a) \to \R$, $c(a) \not= 0$, 
such that $\widetilde{\sigma}(T(t)) = c(t)\sigma(t), (t \in (N, a))$ (see \cite{Mather}). 

\bel
\label{K-equivalence}
The ${\mathcal K}$-equivalence class of a signed area density function $\sigma$ is 
independent of the choice of the frame 
$h_1, h_2, \dots, h_n$ and of the coordinates $t_1, t_2, \dots, t_n$ on $(\R^n, a)$ 
and depend only on the frontal $f$. 
\enl

\Proof
Let us take another frame $k_1, \dots, k_n$. Then 
there exists $A = (a_{ij}) : (\R^n, a) \to {\mbox{\rm GL}}(n, \R)$ 
such that $(h_1, \dots, h_n) = (k_1, \dots, k_n)A$. Then 
$h_1 \wedge h_2 \wedge \cdots \wedge h_n = 
(\det A)(k_1 \wedge k_2 \wedge \cdots \wedge k_n)$. 
Therefore $\sigma$ is transformed to $(\det A)\sigma$. 
Let us take another coordinates $T_1, T_2, \dots, T_n$ on $(\R^n, a)$. 
Then 
$$
(\dfrac{\pa f}{\pa T_1} \wedge \dfrac{\pa f}{\pa T_2} \wedge \cdots \wedge \dfrac{\pa f}{\pa T_n})(T(t)) 
= 
J(t)(\dfrac{\pa f}{\pa t_1} \wedge \dfrac{\pa f}{\pa t_2} \wedge \cdots \wedge \dfrac{\pa f}{\pa t_n})(t), 
$$
where $J(t)$ is the Jacobian function $\pa(T_1, \dots, T_n)/\pa(t_1, \dots, t_n)$ at $t$. 
Therefore $\sigma(t)$ is transformed to the function $J(t)\sigma(T(t))$. 
Thus we have the required result. 
\QED

\

We call the signed density function of a frontal, considered up to ${\mathcal K}$-equivalence, a {\bf density function} of the 
frontal. 
The singular locus (non-immersive locus) $S(f)$ of $f$ coincides with the zero locus $\{ \sigma = 0\}$ 
of the density function $\sigma$.

\section{Proper frontals}
\label{Proper frontals}

Frontals can be collapsing in general. 
For example, any constant mapping $f : N \to M$ is a frontal. In fact any 
lifting $F : N \to \Gr(n, TM)$ of $f$ is Legendre in that case. 
See also Example \ref{global}. 

\bef
{\rm 
A frontal $f : N \to M$ is called a {\bf proper frontal} 
if the regular locus 
$$
R(f) := \{ t \in N \mid f_*: T_tN \to T_{f(t)}M {\mbox{\rm \ is injective.}}\}
$$
of $f$ is dense in $N$. 
A germ of frontal $f : (N, a) \to M$ is called a {\bf germ of proper frontal} 
if there exists a representative of $f$ which is a proper frontal. 
}
\enf

Note that $R(f)$ is an open subset of $N$ in general. Then the condition that 
$f$ is a proper frontal requires that $R(f)$ is open and dense. 

The fundamental property of proper frontals is the following: 

\bep
\label{global-existence}
Let $f : N \to M$ be a proper frontal. Then there exists the unique global Legendre (i.e. $D$-integral) lift 
$\widetilde{f} : N \to \Gr(n, TM)$ of $f$, for 
the canonical projection $\pi : \Gr(n, TM) \to M$, $\pi\circ \widetilde{f} = f$. 
Here $D$ is the contact distribution on $\Gr(n, TM)$, $n = \dim(N)$, 
introduced in \S \ref{Generalised frontals}. 
\enp

\Proof
Consider the mapping $F : R(f) \to \Gr(n, TM)$ defined by 
$F(t) = f_*(T_tN) \in \Gr(n, T_{f(t)}M)$ $\subset \Gr(n, TM)$. Then $F$ is 
a $D$-integral mapping and $\pi\circ F = f\vert_{R(f)}$. 
By Lemma \ref{D-integral}, $F$ is a unique Legendre lifting of $f\vert_{R(f)}$. 
Since $f$ is a frontal, 
for any $a \in N$, there exists an open neighbourhood $U$ of $a$ and 
a $D$-integral lift $\widetilde{f} : U \to \Gr(n, TM)$ of $f$.
Then by the uniqueness of $F$, we have $\widetilde{f} = F$ on $U \cap R(f)$. 
Since $f$ is a proper frontal, $R(f)$ is dense in $N$, and therefore 
$U \cap R(f)$ is dense in $U$. 
Thus the Legendre lift $F$ of $f$ is uniquely extended to $U \cup R(f)$. 
Since $a$ is arbitrary, we have the unique Legendre lift 
$\widetilde{f} : N \to \Gr(n, TM)$ of $f$. 
\QED

\bep
If $f : N \to M$ is a frontal and it has a unique Legendre lift $\widetilde{f} : N \to \Gr(n, TM)$, then $f$ is a proper frontal. 
\enp

\Proof
Suppose the regular locus $R(f)$ of $f$ is not dense in $N$. Then there exists 
a non-void open subset $U \subset N$ such that the maximal rank of $f\vert_U$ is $\ell < n$. 
Then there exists a non-void open subset $V \subset U$ such that $f\vert_V$ is of constant rank $\ell$. 
Then there exists a non-void open subset $W \subset V$ and an open subset $\Omega \subset M$ 
such that $f\vert_W : W \to \Omega$ is right-left equivalent to $h : \R^n \to \R^m$ which is 
defined by $h(s_1, \dots, s_\ell, s_{\ell + 1}, \dots, s_n) = (s_1, \dots, s_\ell, 0, \dots, 0)$ 
(\lq\lq Rank theorem\rq\rq, see \cite{Brocker}). 
Let $\widetilde{f} : N \to \Gr(n, TM)$ be a Legendre lift and $\widetilde{h} : \R^n \to 
\Gr(n, T\R^m)$ be the induced lift of $h$ by $\widetilde{f}\vert_W$ (cf. Proposition \ref{induced-D-integral}). 
Then $T_{h(s)}(\R^\ell\times \{ 0\}) \subset \widetilde{h}(s), (s \in \R^n)$. 
Then there exists a non-trivial perturbation of $\widetilde{h}$ therefore of $\widetilde{f}$ with compact support. 
\QED

\

\ber
\label{non-approximation-}
{\rm
{\it Proper frontals are not generic} in $C^\infty$-topology in general. In fact 
the frontal mapping $f : \R^2 \to \R^3$ constructed in Example \ref{global} can not be approximated by 
any proper frontal. 
}
\enr

Now we introduce the notion of non-degenerate frontals which was originated in 
\cite{KRSUY}. 

\bef
{\rm 
We say that a frontal $f : (N, a) \to M$ has a {\bf non-degenerate} singular point at $a$ if 
the density function $\sigma$ of $f$ satisfies that $\sigma(a) = 0$ and $d\sigma(a) \not= 0$. 
Note that the condition is invariant under the ${\mathcal K}$-equivalence of $\sigma$ (see Proposition \ref{K-equivalence}). 

}
\enf

To study the property of non-degenerate singular points of frontals, we recall 
the following result. 

\bel
\label{determinant}
Let $N$ be a manifold of dimension $n$. 
Let $g : (N, a) \to (N, g(a))$ be a map-germ. Let $J_g$ denote the Jacobi matrix of $g$ and 
$\det(J_g) : (N, a) \to \R$ 
the Jacobian determinant of $g$. Suppose $(\det J_g)(a) = 0$. 
Then $(d\det(J_g))(a) = 0$ if $\rank(J_g)(a) \leq n-2$, that is, if $g$ is of corank $\geq 2$ at $a$. 
\enl

\Proof
It is easy to see, as a fundamental fact in the linear algebra, for the determinant function $\det : M(n, n; \R)$ on the space of $n\times n$-matrices, 
and for any $A \in  M(n, n; \R)$ with $\det(A) = 0$, 
$(d\det)(A) = 0$ if and only if $\rank(A) \leq n-2$. 
Then we have, if $\rank(J_g) \leq n-2$, then $(d\det(J_g))(p) = (J_g)^*(d\det)(p) = 0$. 
\QED

\bel
\label{non-degenerate-proper}
If a frontal $f : (N, a) \to M$ has a non-degenerate singular point at $a$, 
then $f$ is of corank $1$ such that the singular locus 
$S(f) \subset (N, a)$ is a regular hypersurface. 
\enl

\Proof
Let us take a representative $f : U \to M$ of $f$, using the same symbol, satisfying that 
$d\sigma(t) \not= 0$ for any $t \in U$. Then $S(f) = \{ t \in U \mid 
\sigma(t) = 0\}$ is a regular hypersurface of $U$. In particular $S(f)$ is nowhere dense 
in $U$. Therefore $f$ is a proper frontal. 
Let $\widetilde{f} : U \to \Gr(n, TM)$ be the unique Legendre lifting of $f$. 
Set $V = \widetilde{f}(a) \subset T_{f(p)}M$. 
Take a local coordinate system $(x_1, \dots, x_n, x_{n+1}, \dots, x_m)$ 
around $f(a)$ of $M$ such that $V = \langle (\pa/\pa x_1)(f(a)), \dots, (\pa/\pa x_n)(f(a))\rangle_\R$ 
Define $g : U \to \R^n$ by $g = (x_1\circ f, \dots, x_n\circ f)$, deleting $U$ if necessary. 
Then the rank of $g_*$ at $a$ 
is equal to the rank of $f_*$ at $a$. Moreover 
the signed area density function of $g$ is ${\mathcal K}$-equivalent to that of $f$. 
Note that the signed area density function of $g$ is ${\mathcal K}$-equivalent to the 
Jacobian determinant of $g$. 
Suppose the rank of $g_*$ at $p$ is less than $n-1$. 
Then by Lemma \ref{determinant} we see that 
$(d\sigma)(a) = d(\det(J_g))(a) = 0$. This leads a contradiction to the assumption of non-degeneracy. 
Therefore we have that the rank of $g_*$ is equal  to $n-1$. 
Thus we have the required result. 
\QED

\section{Tangent bundles and complementary bundles}
\label{tangent-and-normal}

Let $f : N \to M$ be a proper frontal. 
Let $\dim(N) = n$ and $\widetilde{f} : N \to \Gr(n, TM)$ be the unique Legendre lifting of $f$ (Proposition \ref{global-existence}). 
Then we have a subbundle $T_f$ of the pull-back bundle $f^*TM$ over $N$ defined by 
$$
T_f := \{ (t, v) \in N\times TM \mid v \in \widetilde{f}(t) \} \subset f^*TM. 
$$
We call $T_f$ the {\bf tangent bundle} to the proper frontal $f$. Moreover we call the quotient bundle 
$Q_f := f^*TM/T_f$ the {\bf complementary bundle}. 

\bef
\label{oriented-co-oriented}
{\rm
A proper frontal $f : N \to M$ is called {\bf oriented} (resp. {\bf co-oriented}) if the bundle $T_f$ (resp. $Q_f$) is oriented. 
$f$ is called {\bf orientable} (resp. {\bf co-orientable}) if $T_f$ (resp. $Q_f$) is orientable. 
}
\enf

\bee
{\rm
The proper front $f : S^1 (\subset \C) \to \R^2 (= \C)$ defined by 
$
z \mapsto 2z - \overline{z}^2, 
$
for $z \in \C, \vert z\vert = 1$ (\lq\lq cardioid\rq\rq) is not orientable nor co-orientable. 
The half cuspidal edge $h : \R^2 \setminus \{ (0, 0)\} \to \R^3$ 
(see Example \ref{global}) 
restricted to $\R^2 \setminus \{ (0, 0)\}$ is a proper front which is not orientable nor co-orientable. 
The mapping $\R^2 \to \R^3$ defined by the normal form of the cuspidal edge (resp. folded umbrella) 
is a proper front (resp. frontal) which is orientable and co-orientable. 
}
\ene

Let $f : N \to M$ be a proper frontal. Then 
the bundle homomorphism $\varphi_f : TN \to T_f, \varphi(t, v) = (t, f_*(v))$ is induced. 
Then we have 
$$
R(f) = \{ t \in N \mid f_* : T_tN \to T_{f(t)}M {\mbox{\rm \ is injective}}\} = \{ t \in N \mid 
\varphi_f {\mbox{\rm \ is injective at\ }} t \}. 
$$


\

The notion of frontals will play important role in differential geometry. 
Therefore the following observations are important. 
First we treat the case of hypersurfaces $(m = n + 1)$. 

\bel
\label{frontal-hypersurface}
If $M$ is endowed with an Riemannian metric, then $f : N^n \to M^{n+1}$ is a frontal if and only if, 
for any $a \in N$, there exists an open neighbourhood $U$ of $a$ and a unit vector 
field $\nu$ along $f$ such that $\nu(t)$ is normal to the subspace $f_*(T_tN)$ for any $t \in U$. 
\enl

\Proof
Let $f$ be frontal. Let $\widetilde{f}(t) = (f(t), [\alpha(t)])$ be a Legendre lifting of $f$. 
It defines the local integral tangential hyperplane field $\Ker(\alpha(t))$ along $f$. Then 
we associate the normal line field $\Ker(\alpha(t))^\perp$ 
with $\widetilde{f}(t)$ and take a local unit frame $\nu(t)$ 
of $\Ker(\alpha(t))^\perp$. 
Conversely let $\nu(t)$ be a local unit normal field along $f$ with $f_*(T_tN)$. 
Regarding the metric, we associate a non-zero cotangent vector field $\alpha(t)$ with $\nu(t)$ so 
that $\Ker(\alpha(t)) = \nu(t)^\perp$. 
Then the tangential hyperplane field $\nu(t)^\perp$ satisfies 
the condition $f_*(T_tN) \subseteq \nu(t)^\perp$. The condition is equivalent to that 
$\widetilde{f}(t) = (f(t), [\alpha(t)])$ is a Legendre map. 
\QED

\bel
If $M$ is endowed with an Riemannian metric, then $f : N^n \to M^{n+1}$ is a front if and only if 
locally there exists 
a normal unit vector field $\nu$ along $f$ such that $(f, \nu)$ is an immersion to 
the unit tangent bundle $T_1M$. 
\enl

\Proof
Regard each unit vector $\nu \in T_xM$ as an element of $T_x^*M$ by 
$v \mapsto \nu\cdot v$, we have the natural 
double covering $T_1M \to P(T^*M)$. Therefore we have required result by 
Lemma \ref{frontal-hypersurface}. 
\QED

In generalised cases, we have: 

\bel
If $M$ is endowed with an Riemannian metric, then $f : N^n \to M^m$ is a frontal if and only if, 
for any $a \in N$, there exists an open neighbourhood $U$ of $a$ and a system 
of orthonormal vector fields $\nu_1, \dots, \nu_{m-n}$ over $U$ 
along $f$ such that $\nu_i(t)$ is normal to the subspace $f_*(T_tN)$ for any $t \in U$, $i = 1, \dots, m-n$. 
\enl

\Proof
Suppose $f$ is a frontal. 
For any $a$, let $\widetilde{f} : U \to \Gr(n, TM)$ be a Legendre local lifting of $f\vert_U$. 
Deleting $U$ if necessary, 
take an orthonormal frame $h_1, \dot, h_n, \nu_1, \dots, \nu_{m-n}$ on $U$ 
such that $h_1(t), \dot, h_n(t)$ form a basis of $\widetilde{f}(t) \subset T_{f(t)}M$ 
for any $t \in U$. Then $\nu_1, \dots, \nu_{m-n}$ satisfy the required condition. 
Conversely we may set $\widetilde{f}(t) = \langle \nu_1(t), \dots, \nu_{m-n}(t)\rangle^\perp$. 
Then $\pi\circ\widetilde{f} = f$ and $(\pi\circ\widetilde{f})_*(T_tN) = f_*(T_tN) \subset \widetilde{f}(t)$, 
hence $\widetilde{f}$ is Legendre by Lemma \ref{D-integral}. 
\QED

The following is clear: 

\bel
If $M$ is a Riemannian manifold, then the condition that $f : N^n \to M^m$ is a front 
is equivalent to the local existence of an orthonormal unit frame $\nu_1, \dots, \nu_n$ along $f$ such that 
$t \mapsto (f(t), \langle \nu_1(t), \dots, \nu_n(t)\rangle^\perp)$ is an immersion to 
$\Gr(n, TM)$. 
\enl

\

Let $f : N \to M$ be a proper frontal. 
If $M$ is endowed with a Riemannian metric, then we define the {\bf normal bundle} to $f$ by
$$
N_f := \{ (t, w) \in N\times TM \mid w \in \widetilde{f}(t)^\perp \} \subset f^*TM, 
$$
which is isomorphic to the complementary bundle $Q_f$ (see \S \ref{tangent-and-normal}). 
Note that both bundles $T_f$ and $N_f$ have induced Riemannian bundle structures 
from $TM$.

\section{Openings and frontals}
\label{Openings and frontals} 

In this section, we review the known results on \lq\lq geometric\rq\rq\, openings. 

We denote by ${\mathcal E}_{N, a}$ the $\R$-algebra of $C^\infty$ function-germs on $(N, a)$ 
with the maximal ideal ${\mathfrak m}_{N, a}$. If $(N, a) = (\R^n, 0)$ is the origin, 
then we use ${\mathcal E}_n, {\mathfrak m}_n$ 
instead of ${\mathcal E}_{N, a}, {\mathfrak m}_{N, a}$ respectively. 

\bef
\label{Jacobi-module}
{\rm 
(\cite{Ishikawa92}\cite{Ishikawa96}) 
Let $f : (N, a) \to (M, b)$ be a $C^\infty$ map-germ with $\dim(N) = n \leq m = \dim(M)$. 
We define the {\bf Jacobi module} of $f$: 
$$
{\mathcal J}_f  :=  {\mathcal E}_{N, a}\, d(f^*\Omega_{M, b}) 
 = \{ \ \sum_{j=1}^m \ a_j \ df^j \mid \ a_j \in {\mathcal E}_{N, a}, 1 \leq j \leq m \} 
$$
in the space $\Omega_{N, a}^1$ of $1$-form germs on $(N, a)$. Here $f^j = x^j\circ f$, for a system of coordinates 
$(x^1, \dots, x^m)$ of $(M, b)$. 
Further we define the {\bf ramification module} ${\mathcal R}_f$ by 
$$
{\mathcal R}_f := \{ h \in {\mathcal E}_{N, a} \mid dh \in {\mathcal J}_f \}. 
$$
}
\enf

\bee
\label{example-1-variable}
{\rm 
Let $\mu$ be a positive integer and $g : (\R, 0) \to (\R, 0), g(t) = t^\mu$. 
Then ${\mathcal J}_g = {\mathfrak m}_1^{\mu - 1}dt$ and 
${\mathcal R}_g = \R + {\mathfrak m}_1^\mu$. 
Here ${\mathfrak m}_1^\mu = t^\mu{\mathcal E}_1 = \{ h \in {\mathcal E}_1 \mid \frac{d^{k}h}{dt^k}(0) = 0, (0 \leq k \leq \mu) \}$. 
In fact, since $dg = \mu t^{\mu - 1}dt$, we gave ${\mathcal J}_g = {\mathfrak m}_1^{\mu - 1}dt$. 
Moreover, for a $k \in {\mathcal E}_1$, we have that $k \in {\mathcal R}_f$ if and only if $\frac{dk}{dt} \in 
{\mathfrak m}_1^{\mu - 1}$ if and only if $k \in \R + {\mathfrak m}_1^\mu$. 
}
\ene

Note that ${\mathcal J}_f$ is just the first order component of the graded differential ideal 
${\mathcal J}^{\bullet}_f$ in $\Omega_{N, a}^{\bullet}$ 
generated by $df^1, \dots, df^m$. 
Then the singular locus is given by 
$
S(f) = \{ x \in (N, a) \mid \rank\,{\mathcal J}_f(x) < n \}. 
$
Also we consider the {\bf kernel field} 
$
\Ker(f_* : TN \to TM), 
$
of $f$ near $a$.  
Then we see that, for another map-germ $f' : (N, a) \to (M', b')$ with 
${\mathcal J}_{f'} = {\mathcal J}_{f}$, $n \leq m'$, we have 
$
\Sigma_{f'} = \Sigma_{f}
$
and $\Ker(f'_*) = \Ker(f_*)$. 
Note that related notion was introduced in \cite{Mond3}. 

\bel
Let $f : (N, a) \to (M, b)$ be a map-germ. Then we have:  
\\
{\rm (1)} 
 $f^*{\mathcal E}_{M,b} \subset {\mathcal R}_f \subset {\mathcal E}_{N,a}$ and 
${\mathcal R}_f$ is an ${\mathcal E}_{M,b}$-module via $f^*$. 
\\
{\rm (2)}
For another map-germ $f' : (N, a) \to (M', b')$, 
${\mathcal J}_{f'} = {\mathcal J}_f$ if and only if ${\mathcal R}_{f'} = {\mathcal R}_f$. 
\\
{\rm (3)}
If $\tau : (M, b) \to (M', b')$ is a diffeomorphism-germ, then 
${\mathcal R}_{\tau\circ f} = {\mathcal R}_f$. If $\sigma : 
(N', a') \to (N, a)$ is a diffeomorphism-germ, then 
${\mathcal R}_{f\circ\sigma} = \sigma^*({\mathcal R}_f)$. 
\enl

\Proof (1) follows from that,  
if $h \in {\mathcal R}_f$ and $dh = \sum_{j=1}^m p_j df_j$, then we have
$$
d\{(k\circ f)h\} = \sum_{j=1}^m\left\{ (k\circ f)p_j + h\left({\pa k}/{\pa y_j}\right)\right\} df_j. 
$$

(2) It is clear that ${\mathcal J}_{f'} = {\mathcal J}_f$ implies 
${\mathcal R}_{f'} = {\mathcal R}_f$. Conversely suppose ${\mathcal R}_{f'} = {\mathcal R}_f$. 
Then any component $f'_j$ of $f'$ belongs to 
${\mathcal R}_{f'} = {\mathcal R}_f$, hence $df_j \in {\mathcal J}_f$. 
Therefore ${\mathcal J}_{f'} \subset {\mathcal J}_f$. 
By the symmetry we have ${\mathcal J}_{f'} = {\mathcal J}_f$. 

(3) follows from that 
${\mathcal J}_{\tau\circ f} = {\mathcal J}_f$ and 
${\mathcal J}_{f\circ\sigma} = \sigma^*({\mathcal J}_f)$. 
\QED

\bef
\label{J-equivalence}
{\rm
Let $f : (N, a) \to M$ and $g : (N',a') \to M'$ be map-germs. 
Then $f$ and $g$ are called {\bf ${\mathcal J}$-equivalent} if 
there exists a diffeomorphism-germ $\sigma : (N, a) \to (N', a')$ such that 
${\mathcal J}_{g\circ\sigma} = {\mathcal J}_f$. 
Here ${\mathcal J}_f = {\mathcal E}_{N, t_0}f^*\Omega^1_{M, f(t_0)}$ (see Definition \ref{Jacobi-module}). 
Note that $\dim(M)$ and $\dim(M')$ can be different. 
}
\enf

\bef
{\rm 
Let $f : (\R^n, a) \to (\R^m, b)$ a map-germ and 
$h_1, \dots, h_r \in {\mathcal R}_f$. Then the map-germ 
$F : (\R^n, a) \to \R^m\times\R^r = \R^{m+r}$ defined by 
$$
F = (f_1, \dots, f_m, h_1, \dots, h_r)
$$
is called an {\bf opening} of $f$, while $f$ is called a {\bf closing} of $F$. 
}
\enf

\bep
Let $M, N$ be a manifold of dimension $m, n$ respectively. 
A map-germ $f : (N, a) \to M$ is a frontal if and only if 
$f$ is right-left equivalent to an opening of a map-germ $g : (\R^n, 0) \to (\R^n, 0)$. 
\enp

\Proof
Suppose  $f : (N, a) \to M$ be a frontal map-germ. 
Then, since $\widetilde{f}(a)$ is an $n$-dimensional vector subspace of $T_{f(a)}M$, 
there exists of a system of local coordinates of $(M, f(a))$ 
$$
y_1, \dots, y_n, z_1, \dots, z_k, (k = m - n), 
$$
such that, for  
$g = (y_1\circ f, \dots, y_n\circ f)$, the component $z_j\circ f, (1 \leq j \leq k)$ belongs to the 
${\mathcal R}_g$ and therefore $f$ is right-left equivalent to an opening of $g$. 
Conversely, suppose $f$ is right-left equivalent to an opening $G : (\R^n, a) \to \R^{n + k} = \R^m$ of a germ $g = (g_1, \dots, g_n) : (\R^n, a) \to (\R^n, g(a))$. 
Set $G = (g_1, \dots, g_\ell, h_1, \dots, h_k)$. Then, since $h_j \in {\mathcal R}_g, 1 \leq j \leq n$
$dh_j = \sum_{i=1}^n a_j^idg_i$, for some function-germs $a_j^i : (\R^n, a) \to \R$. 
Define $\widetilde{G} : (\R^n, a) \to \Gr(n, T\R^m)$, in terms of Grassmannian coordinates, 
$$
\widetilde{G}(t) = \left(g_1(t), \dots, g_\ell(t), h_1(t), \dots, h_k(t), a_j^i(t)\right), \quad(t \in (\R^n, a)). 
$$
Then, by Lemma \ref{canonical-coordinate}, 
$\widetilde{G}$ is a $D$-integral lift of $g$. Therefore $G$ is an $\ell$-frontal, and so is $f$. 
\QED

\section{Versal openings}
\label{Versal openings}

\bef
\label{versal-opening-def}
{\rm 
An opening $F = (f, h_1, \dots, h_r)$ of $f$ is called a {\it versal opening} (resp. 
a {\it mini-versal opening}) of $f : (\R^n, a) \to (\R^m, b)$, if 
$1, h_1, \dots, h_r$ form a (minimal) system of 
generators of ${\mathcal R}_f$ as an ${\mathcal E}_{\R^m, b}$-module via 
$f^* : {\mathcal E}_{\R^m, b} \to {\mathcal E}_{\R^n, a}$. 
}
\enf

A $C^\infty$ map-germ $f : (\R^n, a) \to (\R^m, b)$ 
is called {\bf analytic} if $f$ is right-left equivalent to a real analytic map-germ (\cite{Ishikawa96}). 
Moreover $f$ is called a {\bf finite map-germ} if 
${\mathcal E}_{\R^n, a}$ is a finite $f^*({\mathcal E}_{\R^m, b})$-module. 
Then $f$ is finite if and only if 
$\dim_{\R}({\mathcal E}_{\R^n, a}/\langle f_1, \dots, f_m\rangle_{{\mathcal E}_{\R^n, a}} < \infty$. 
If $f$ is analytic, then $f$ is finite if and only if its complexification has isolated zero set (\cite{Wall81}). 
By ${\mathfrak m}_{\R^m, a}$, we denote the maximal ideal of ${\mathcal E}_{\R^m, a}$ 
which consists of function-germs vanishing at $0$. 
By the projection $\pi_m : \R^{m+r} = \R^m\times\R^r \to \R^m$ we regard 
$\R^{m+r}$ as an affine bundle over $\R^m$. 
If $f$ is finite and analytic, then, in the analytic category, 
${\mathcal R}^\omega f$ is a finite ${\mathcal O}_{\R^m, b}$-module. 

We summarise the known results on the existence of versal openings: 


\bet
\label{main theorem}
Let $f : (\R^n, a) \to (\R^m, b)$ be a map-germ. 
Suppose that {\rm (I)} $f$ is finite and of corank at most one, or 
{\rm (II)} $f$ is finite analytic. 
Then we have
\\
{\rm (1)} The ramification module ${\mathcal R}_f$ of $f$ is a finitely generated 
$f^*({\mathcal E}_{\R^m, b})$-module. In particular $f$ has a versal opening. 
\\
{\rm (2)} 
$1, h_1, \dots, h_r \in {\mathcal R}_f$ form a system of generators of ${\mathcal R}_f$ as a 
$f^*({\mathcal E}_{\R^m, b})$-module if and only if 
the residue classes 
$1, \overline{h}_1, \dots, \overline{h}_r$ 
form an $\R$-basis 
of the vector space 
$V = {\mathcal R}_f/f^*({\mathfrak m}_{\R^m, b}){\mathcal R}_f$. 
In particular there exists a versal opening $F = (f, h_1, \dots, h_r)$ of $f$. 
If $r = \dim_{\R}V - 1$, then $F$ is called a {\it mini-versal opening} of $f$. 
\\
{\rm (3)} For any versal opening $F : (\R^n, a) \to \R^{m+r}$ of $f$ and for any 
opening $G : (\R^n, a) \to \R^{m+s}$ of $f$, 
there exists an affine bundle map $\Psi : (\R^{m+r}, F(a)) \to (\R^{m+s}, G(a))$ 
such that $G = \Psi\circ F$. 
\\
{\rm (4)} For any two mini-versal opening $F, F' : (\R^n, a) \to \R^{m+r}$ of $f$, 
there exists an affine bundle isomorphism $\Phi : (\R^{m+r}, F(0)) \to (\R^{m+r}, F'(0))$ such that 
$F' = \Phi\circ F$. 
\ent

\noindent
{\it Proof of Theorem \ref{main theorem}:} 
Case (I): \ (1) is proved as Lemma 2.1 of \cite{Ishikawa96}. 
Then $f^* : {\mathcal E}_{\R^m, b} \to {\mathcal R}_f$ is 
a homomorphism of differentiable algebras in the sense of Malgrange \cite{Malgrange}. 
Therefore, by Malgrange's preparation theorem (\cite{Malgrange} Corollary 4.4) we have (2). 
Case (II): Let $1, h_1, \dots, h_r$ generate 
${\mathcal R}^\omega f$ over ${\mathcal O}_{\R^m, b}$ via $f^*$. 
Then $1, h_1, \dots, h_r$ generate ${\mathcal R}_f$ over ${\mathcal E}_{\R^m, b}$ via $f^*$ 
by Proposition 5.2 of \cite{Ishikawa14}. Therefore we have (1) and (2). 
The assertion (3) is clear from the definitions. (4) follows from (2).  \QED

\

We do not repeat the proofs of Lemma 2.1 in \cite{Ishikawa96} nor 
Proposition 5.2 in \cite{Ishikawa14}. However we give an exposition, relating the theory of $C^\infty$-rings, on 
Malgrange's preparation theorem in 
\S \ref{Appendix: Malgrange preparation theorem on differentiable algebras} of this paper. 

%
%
%
%
%


\section{Subfrontals and superfrontals}
\label{Subfrontals and superfrontals}

Relating the theory of frontals with that of openings, we are led to 
the following generalisations of frontals naturally. 

\bef
{\rm 
Let $M$ be a manifold and $\ell$ an integer with $0 \leq \ell \leq \dim(M)$. 
A map-germ $f : (N, a) \to M$ is called an {\bf $\ell$-frontal} 
if there exists a $D$-integral lift $\widetilde{f} : (N, a) \to \Gr(\ell, TM)$ 
of $f$. 
Here we do not assume that $\ell = \dim(N)$ and $D$ is the contact distribution 
on $\Gr(\ell, TM)$. 
The condition on the $C^\infty$ mapping $\widetilde{f}$ is that 
$f_*(T_tN) \subset \widetilde{f} \in \Gr(\ell, T_{f(t)}M)$ for any $t \in (N, a)$. 
If $0 < \ell < \dim(N)$, then an $\ell$-frontal is called a {\bf subfrontal}. 
If $\dim(N) < \ell < \dim (M)$, then an an $\ell$-frontal is called a {\bf superfrontal}. 
}
\enf

\bep
Let $f : (N, a) \to M$ be an $\ell$-frontal and $f' : (N', a') \to M'$ be right-left equivalent to $f$. 
Then also $f'$ is an $\ell$-frontal. 
\enp

\Proof
The proof is performed similarly to Proposition \ref{induced-D-integral} using Proposition \ref{canonical-is-canonical}. 
\QED

\bep
Let $M, N$ be a manifold of dimension $m, n$ respectively and $\ell$ an integer with $0 \leq \ell \leq m$. 
A map-germ $f : (N, a) \to M$ is an $\ell$-frontal if and only if 
$f$ is right-left equivalent to an opening of a map-germ $g : (\R^n, 0) \to (\R^\ell, 0)$. 
\enp

\Proof
Suppose  $f : (N, a) \to M$ be an $\ell$-frontal map-germ. 
Then, since $\widetilde{f}(a)$ is an $\ell$-dimensional vector subspace of $T_{f(a)}M$, 
there exists of a system of local coordinates of $(M, f(a))$ 
$$
y_1, \dots, y_\ell, z_1, \dots, z_k, (k = m - \ell), 
$$
such that, for  
$g = (y_1\circ f, \dots, y_\ell\circ f)$, the component $z_j\circ f, (1 \leq j \leq k)$ belongs to the 
${\mathcal R}_g$ and therefore $f$ is right-left equivalent to an opening of $g$. 
Conversely, suppose $f$ is right-left equivalent to an opening $G : (\R^n, a) \to \R^{\ell + k} = \R^m$ of a germ 
$g = (g_1, \dots, g_\ell) : (\R^n, a) \to (\R^\ell, g(a))$. 
Set $G = (g_1, \dots, g_\ell, h_1, \dots, h_k)$. Then, since $h_j \in {\mathcal R}_g, 1 \leq j \leq \ell$
$dh_j = \sum_{i=1}^\ell a_j^idg_i$, for some function-germs $a_j^i : (\R^n, a) \to \R$. 
Define $\widetilde{G} : (\R^n, a) \to \Gr(\ell, T\R^m)$, in terms of Grassmannian coordinates, 
$$
\widetilde{G}(t) = \left(g_1(t), \dots, g_\ell(t), h_1(t), \dots, h_k(t), a_j^i(t)\right), \quad(t \in (\R^n, a)). 
$$
Then, by Lemma \ref{canonical-coordinate}, 
$\widetilde{G}$ is a $D$-integral lift of $g$. Therefore $G$ is an $\ell$-frontal, and so is $f$. 
\QED

\section{Algebraic openings}
\label{Algebraic openings}

In this section we will utilise the notion of sheaves which describes locally defined objects (\cite{Bredon}) 
to introduce an algebraic notion which is related to frontals. 

Let $N$ be a manifold. Let ${\mathcal E}_N$ denote the sheaf of $C^\infty$ function-germs on 
$N$. For any open subset $U \subset N$, ${\mathcal E}_N(U) = C^\infty(U)$, 
the $C^\infty$-ring of all real-valued $C^\infty$ functions on $U$. 
Note that ${\mathcal E}_N$ has the natural structure of $C^\infty$-ring sheaf. 
See \S \ref{Appendix: Malgrange preparation theorem on differentiable algebras} for the notion of 
$C^\infty$-rings. 

\bef
{\rm 
Let ${\mathcal F}$ be a sub $C^\infty$-ring sheaves of ${\mathcal E}_N$. 
The {\bf versal opening} $\widetilde{\mathcal F}$ of ${\mathcal F}$ is defined as follows: 
For any open subset $U \subset N$, $\widetilde{\mathcal F}(U)$ is the 
set of $h \in {\mathcal E}_N(U)$ satisfying that, 
for any $p \in U$, there exists 
$g_1, \dots, g_r \in {\mathcal F}_p$ and $a_1, \dots, a_r \in {\mathcal E}_{N,p}$ such that 
$$
dh = \sum_{i=1}^r a_i dg_i 
$$
in $\Omega_{N,p}^1$. Here $\Omega_N^1$ means the sheaf of $C^\infty$ $1$-form-germs 
on $N$, ${\mathcal E}_{N,p}$ (resp. $\Omega_{N,p}^1$) the stalk of ${\mathcal E}_N$ (resp. of $\Omega_N^1$) at $p$, 
i.e. the set of germs at $p$, 
and 
$d : {\mathcal E}_{N,p} \to \Omega_{N,p}^1$ the exterior differential. 

Let ${\mathcal F}, {\mathcal G}$ be a sub $C^\infty$-ring sheaves of ${\mathcal E}_N$. 
Then ${\mathcal G}$ is called an {\bf opening} of ${\mathcal F}$ if 
${\mathcal F} \subseteq {\mathcal G} \subseteq \widetilde{\mathcal F}$. 
}
\enf

We have that following basic properties of algebraic openings. 

\bep
\label{algebraic-opening}
Let $N$ be a $C^\infty$ manifold and let ${\mathcal E}_N$ denote the sheaf of $C^\infty$ function-germs on $N$. 
Let ${\mathcal F}$ be a sub $C^\infty$-ring sheaf of ${\mathcal E}_N$. Then we have
\\
{\rm (1)} 
$\widetilde{\mathcal F}$ is a sub $C^\infty$-ring sheaf of ${\mathcal E}_N$. \ 
{\rm (2)} ${\mathcal F} \subset \widetilde{\mathcal F}$. \ 
{\rm (3)}
$\widetilde{\hspace{-0.1truecm}\widetilde{\mathcal F}} = \widetilde{\mathcal F}$. 
\enp

\Proof
(1) 
Let $p \in N$. 
Let $h_1, \dots, h_r \in \widetilde{\mathcal F}_p$ and $f \in C^\infty(\R^r)$. 
Let $dh_i = \sum_{j=1}^{s_i} a_{ij}dg_{ij}$ for some $a_{ij} \in {\mathcal E}_{N, p}, \, g_{ij} \in {\mathcal F}_p$. 
Then 
$$
d(f(h_1, \dots, h_r)) = \sum_{i=1}^r \frac{\pa f}{\pa x_i}(h_1, \dots, h_r)\, dh_i 
= \sum_{i=1}^r \sum_{j=1}^{s_i} (\frac{\pa f}{\pa x_i}(h_1, \dots, h_r)a_{ij})dg_{ij}. 
$$
Therefore $f(h_1, \dots, h_r)  \in \widetilde{\mathcal F}_p$. 
(2) Let $p \in N$ and $g \in {\mathcal F}_p$. Then we have 
$dg = 1\cdot dg$, and therefore $g \in \widetilde{\mathcal F}_p$. 
(3) 
Let $p \in N$ and $h \in \widetilde{\hspace{-0.1truecm}\widetilde{\mathcal F}}_p$. 
Then $dh = \sum_{i=1}^r a_i dh_i$ for some $a_i \in {\mathcal E}_{N, p}, h_i \in \widetilde{\mathcal F}_p$. 
Since $h_i \in \widetilde{\mathcal F}_p$ for each $i$, $dh_i = \sum_{j=1}^{s_i} b_{ij} dg_{ij}$ for some 
$b_{ij} \in {\mathcal E}_{N, p}$ and $g_{ij} \in {\mathcal F}_p$. Then we have 
$dh = \sum_{i=1}^r \sum_{j=1}^{s_i} (a_i b_{ij}) dg_{ij}$, therefore $h \in \widetilde{\mathcal F}_p$. 
\QED

\

We call ${\mathcal F}$ {\bf full} if $\widetilde{\mathcal F} = {\mathcal F}$. 
Then Proposition \ref{algebraic-opening} shows that 
$\widetilde{\mathcal F}$ is the minimal full sheaf containing ${\mathcal F}$. 

\

Let $\varphi : N' \to N$ be a $C^\infty$ mapping 
and ${\mathcal F}$ a subsheaf of ${\mathcal E}_N$ on $N$. 
We define a subsheaf $\varphi^*{\mathcal F}$ of ${\mathcal E}_{N'}$ on $N'$
by $(\varphi^*{\mathcal F})_q = \varphi^*({\mathcal F}_{\varphi(q)})$, 
where $\varphi^* : {\mathcal E}_{N,\varphi(q)} \to {\mathcal E}_{N',q}$ is defined by 
$\varphi^*(h) = h\circ\varphi, (h \in {\mathcal E}_{N,\varphi(q)})$. 
If $\varphi = \Phi$ is a diffeomorphism, then 
$(\Phi^*{\mathcal F})(U') = \Phi^*({\mathcal F}(\Phi(U')))$ for any open $U' \subset N'$. 

Then we have the naturality of versal openings: 

\bep
\label{algebraic-opening2}
Let ${\mathcal F}$ be a sub $C^\infty$-ring sheaves of ${\mathcal E}_N$. 
For any diffeomorphism $\Phi : N' \to N$ from another manifold $N'$, we have
$\widetilde{\Phi^*{\mathcal F}} = \Phi^*{\widetilde{\mathcal F}}$. 
\enp

\Proof
Let $q \in N'$ and $h \in \widetilde{\Phi^*{\mathcal F}}_q$. 
Then $dh = \sum_{i=1}^r a_i d(\Phi^*g_i)$ for some $a_i \in {\mathcal E}_{N', q}$ and $g_i \in {\mathcal F}_{\Phi(q)}$. 
Then $d(\Phi^{-1*}h) = \Phi^{-1*}(\sum_{i=1}^r a_i d(\Phi^*g_i)) = \sum_{i=1}^r (\Phi^{-1*} a_i) dg_i$. 
Since $\Phi^{-1*} a_i \in {\mathcal E}_{N, \Phi(q)}$, we see $\Phi^{-1*}h \in {\widetilde{\mathcal F}}_{\Phi(q)}$, 
therefore $h \in \Phi^*(({\widetilde{\mathcal F}})_{\Phi(q)}) = (\Phi^*{\widetilde{\mathcal F}})_q$. 
Thus we have $\widetilde{\Phi^*{\mathcal F}}_q \subseteq (\Phi^*{\widetilde{\mathcal F}})_q$. Applying the 
same argument to $\Phi^{-1}$ and $\Phi^*{\mathcal F}$, then we have $\widetilde{\Phi^*{\mathcal F}}_q \supseteq (\Phi^*{\widetilde{\mathcal F}})_q$. Therefore we have the required equality. 
\QED

\

\bef
{\rm
Let $N$ be a manifold. Let ${\mathcal F}$ be a sub $C^\infty$-ring sheaf of ${\mathcal E}_N$. 
A mapping $f : N \to M$ is called a {\bf realisation} of ${\mathcal F}$ if 
${\mathcal F} = f^*{\mathcal E}_M$. 
}
\enf

The following is clear: 

\bep
Let $f : (\R^n, a) \to (\R^m, b)$ be a  map-germ. Let ${\mathcal F} = f^*{\mathcal E}_{\R^m, b}$ 
be the germ of subsheaf of ${\mathcal E}_{\R^n, a}$. 
Let $F : (\R^n, a) \to \R^{m+r}$ be an opening of $f$. 
Then $F$ is a versal opening of $f$ if and only if $F$ is a realisation of the algebraic opening $\widetilde{\mathcal F}$ 
of ${\mathcal F}$. 
\enp

%
%

\bef
{\rm 
A mapping $f : N \to M$ is called {\bf locally injective} 
if for any $a \in N$, there exists an open neighbourhood $U$ of $a$ in $N$ such that 
$f\vert_U : U \to M$ is injective. 
}
\enf

\bep
Let $f : N \to M$ be a finite mapping and $F : N \to M'$ a realisation of 
the versal opening $\widetilde{f^*{\mathcal E}_M}$ of $f^*{\mathcal E}_M$. 
Then $F$ is locally injective. 
\enp

\Proof
Let $a \in N$. Then $F^*{\mathcal E}_{M', F(a)} = \widetilde{f^*{\mathcal E}_M}_a = {\mathcal R}_{f, a}$. 
Then the germ of $F$ at $a$ is a versal opening of $f$. Therefore by Proposition 2.16 of \cite{Ishikawa14} 
we have the result. 
\QED

\bef
{\rm 
Let ${\mathcal F}$ be a sub $C^\infty$-ring of ${\mathcal E}_N$. 
We call ${\mathcal F}$ {\bf locally injective} if 
for any $a \in N$, there exist $h_1, \dots, h_r \in {\mathcal F}_a$
such that $(h_1, \dots, h_r) : (N, a) \to \R^r$ has an injective representative. 
}
\enf

\bep
If $f : N \to M$ is a realisation of a locally injective 
sub $C^\infty$-ring ${\mathcal F}$ of ${\mathcal E}_N$, then $f$ is locally injective. 
\enp

\Proof
Let $a \in N$. There exist $h_1, \dots, h_r \in {\mathcal F}_a$
such that $(h_1, \dots, h_r) : (N, a) \to \R^r$ has an injective representative. 
There exists a $g_i \in {\mathcal E}_{M, f(a)}$ such that $h_i = g_i\circ f$ for each $i, 1 \leq i \leq r$. 
After taking representatives of germs we have 
$h_i = g_i\circ f : U \to \R, (1 \leq i \leq r)$ on an open neighbourhood of $a$. 
Deleting $U$ if necessary, $(h_1, \dots, h_r) = (g_1, \dots, g_r)\circ f : U \to \R^r$ is injective. Therefore 
$f\vert_U$ is injective. 
\QED

\

\

\begin{center}
{\bf {\LARGE Part II.  Advanced studies and applications}}
\end{center}

\section{Frontal curves}
\label{Frontal curves}

Let us give several observations on frontal map-germs and frontal maps $N \to M$ with $\dim(N) = 1$. 

Let $f : (N, a) \to M$ be a map-germ with $\dim(N) = 1$. 
We consider the classification problem of germs up to the right-left equivalence. 
To simplify this, let $(N, a) = (\R, 0)$ and $(M, f(a)) = (\R^m, 0)$. 
Let $t$ be the coordinate of $(\R, 0)$ and $x_1, \dots, x_m$ of $(\R^m, 0)$. 
We define the {\bf order} of $f$ at $0$ by 
$$
\ord(f) := \inf\left\{ k \in \NN \ \left\vert \ \frac{d^{k}f}{dt^k}(0) \not= 0 \right\}\right.
$$
If the Taylor infinite series of $f$ is $0$, then we set $\ord(f) = \infty$. 
It is easy to see that $\ord(f)$ is invariant under right-left equivalence. 

\bel
If $\ord(f) < \infty$, then $f$ is a frontal. 
Moreover $f$ is right-left equivalent to an opening of the map-germ $g : (\R, 0) \to (\R, 0)$ 
defined by $t \mapsto t^\mu$, where $\mu = \ord(f)$. 
\enl

\Proof
For a diffeomorphism-germ $\sigma : (\R, 0) \to (\R, 0)$ and 
a linear transformation $\Phi : (\R^m, 0) \to (\R^m, 0)$, $\Phi\circ f\circ\sigma$ is of form:
$$
\Phi\circ f\circ\sigma = (t^\mu, \ h_2(t), \ \dots\ , \ h_m(t)), 
$$
with $h_i \in {\mathfrak m}_1^{\mu + 1}, 2 \leq i \leq m$. 
Set $g(t) = t^\mu$. Then ${\mathcal R}_g = \R + {\mathfrak m}_1^\mu$ and $h_i \in {\mathcal R}_g, 2 \leq i \leq m$ 
(see Example \ref{example-1-variable}). 
Therefore $\Phi\circ f\circ\sigma$ is an opening of $g$. 
\QED

\bec
Let 
$f : (\R, a) \to \R^m$ be an analytic map-germ. If $f$ is not a constant map-germ, then $f$ is a frontal. 
\enc

As for a global result, we have:

\bel
\label{global-n=1}
Let $\dim(N) = 1$ and $f : N \to M$ a frontal. 
Then there exists a global Legendre lift $\widetilde{f} : 
N \to M$. 
\enl

\Proof
Let $R(f)$ denote the immersion locus of $f$ and set $S := N \setminus \overline{R(f)}$. 
We have the Legendre lift 
$F : R(f) \to \Gr(1, TM)$ of $f\vert_{R(f)}$ which is defined by $F(t) = f_*(T_tN)$. 
The mapping $F$ is extended to $\overline{R(f)}$ continuously. 
Since $f$ is a frontal. $F$ is extended to a $C^\infty$ Legendre lift of $f$ 
on an open neighbourhood of $\overline{R(f)}$. 
Now take any connected component $J$ of the open set $S$. 
Then $J$ is diffeomorphic to $S^1$ or an open interval. 
In the case that $J$ is diffeomorphic to $S^1$, then $f\vert_J$ is of constant rank $0$ and 
it is a constant mapping. Let us consider the case that 
$J \subset N$ is diffeomorphic to an open interval. Take the closure $I = \overline{J}$ in $N$, 
which is diffeomorphic to an interval, $[0, 1], (0, 1]$ or $(0, 1)$. 
In the case $I$ is diffeomorphic to $[0, 1]$, consider the boundary points of $I$, which 
belong to $R(f)$ necessarily. Since the fibre of $\pi : \Gr(1, TM) \to M$ is 
diffeomorphic to the projective space $\Gr(1, \R^m) = \R P^{m-1}$ is connected, 
we can extend the given Legendre lift $F : R(f) \to \Gr(1, TM)$ to 
a Legendre lift an open set containing $\overline{R(f)} \cup I$. 
The extension is performed 
independently for each connected component of $S$. 
Thus we have a global Legendre lift $\widetilde{f} : 
N \to M$. 
\QED

\ber
\label{non-global}
{\rm 
If $\dim(N) =2$, then a frontal $f : N \to M$ need not have a global Legendre lifting. See Example \ref{global}. 
}
\enr

\

Next we study the genericity problem of frontal curves. To simplify the story we treat 
frontals $f : \R \to \R^m$. Let $\widetilde{f} : \R \to \Gr(1, T\R^m) = P(T\R^m) = \R^m\times \R P^{m-1}$ 
be an integral lifting of $f$ (see Lemma \ref{global-n=1}). Then, turning upside-down the view point, we start from 
an integral map $F : \R \to \Gr(1, T\R^m)$. Let $F$ be, in terms of Grassmannian coordinates $x^1, \dots, x^m, a^2, \dots, a^m$, 
$$
F(t) = (x^1(t), \dots, x^m(t), a^2(t), \dots, a^m(t)), 
$$
which satisfies $F^*\theta^2 = 0, \dots, F^*\theta^m = 0$, namely that
$$
dx^2 - a^2dx^1 = 0, \dots, dx^m - a^mdx^1 = 0. 
$$
The condition is equivalent to that 
$$
\frac{d x^2}{dt}(t) = a^2(t)\frac{d x^1}{dt}(t), \ \dots, \ \frac{d x^m}{dt}(t) = a^m(t)\frac{d x^1}{dt}(t). 
$$
Therefore, if functions $x^1(t), a^2(t), \dots, a^m(t)$ and values $x^2(0), \dots, x^m(0)$ are arbitrarily given, 
then the integral mapping $F$ is uniquely determined. 
Thus we can apply ordinary transversality theorem to discuss the genericity of frontal curves through Legendre curves.

\ber
{\rm
In general we can apply transversality argument to Legendre mappings of corank $\leq 1$ and obtain 
the classification of generic singularities (see \cite{Ishikawa94}\cite{Ishikawa05}). However the similar argument does not work 
for Legendre mappings having singularities of corank $\geq 2$ (see Example \ref{global}, Remark \ref{non-approximation-}). 
}
\enr

\section{Frames and flags}
\label{Frames and flags}

As refinements of the notion of frontal curves, we consider framed curves or \lq\lq flagged\rq\rq\, curves. 
Flagged curves and framed curves in a space-form plays important roles in topology, geometry and singularity theory. 
For example, as it is well-known, the self-linking number in 
$3$-space is defined via framing (\cite{Pohl}). 
The fundamental theory of curves is formulated via osculation framing. Surface boundaries have adapted framings, etc. 
Two kinds of frames, adapted frames and osculating frames, are considered in \cite{Ishikawa12-2} 
from the viewpoint of duality. We classify the singularities of envelopes associated to framed curves.  
The singularities of envelopes in $E^3$ were studied in \cite{Ishikawa10} to apply to  
the flat extension problem of a surface with boundary. 
The problem on extensions by tangentially degenerate surfaces motivates to study the envelopes 
associated to framings on curves in a space form. 


In this article already we have used Grassmannians to introduce the frontals. 
Then we are naturally led to the following definitions. 

Let $M$ be a manifold of dimension $m$ and $\ell_1, \dots, \ell_r$ integers with $0 < \ell_1 < \cdots < \ell_r < m$. 
Define the {\bf flag bundle} $\Fl(\ell_1, \dots, \ell_r; TM)$ over $M$ of type $(\ell_1, \dots, \ell_r)$ as 
the totality of flags $V_{\ell_1} \subset \cdots \subset V_{\ell_r} \subset T_xM$ with $\dim(V_{\ell_i}) = \ell_i, (1 \leq i \leq r)$, 
$x$ running over $M$. Then $\pi : \Fl(\ell_1, \dots, \ell_r; TM) \to M$ is a fibration with fibres of dimension
$$
\ell_1(m - \ell_1) + (\ell_2 - \ell_1)(m - \ell_2) + \cdots + (\ell_r - \ell_{r-1})(m - \ell_r)
$$
Moreover $\Fl(n ; TM) = \Gr(n, TM)$. 

Set ${\mathcal F} = \Fl(\ell_1, \dots, \ell_r; TM)$. 
Suppose $M$ is endowed with an affine connection $\nabla$. 
Let $\gamma(t) = (x(t), V_{\ell_1}(t), \dots, V_{\ell_r}(t))$ be a curve on ${\mathcal F}$. 
Let vectors $v_1(t), \dots v_{\ell_r}(t) \in T_{x(t)}M$ satisfy that 
$V_{\ell_j}(t) = \langle v_1(t), \dots, v_{\ell_j}(t)\rangle_{\R}$ for each $t$ and $1 \leq j \leq r$. 
Then consider the condition 
$$
x'(t) \in V_{\ell_1}(t), \  \nabla v_1(t), \dots, \nabla v_{\ell_j}(t) \in V_{\ell_{j+1}}(t), (1 \leq j < r), 
$$
at $t$. Here, for a vector field $v(t)$ along a curve $x(t)$ in $M$, we define $\nabla v(t) := \nabla_{x'(t)} v(t)$, 
the covariant derivative of $v(t)$ by the velocity vector $x'(t)$. 
By this condition we define the distribution $D \subset T{\mathcal F}$, which depends on the given affine connection. 

If $M$ is a projective space, then the above construction is more clarified (\cite{Ishikawa12-2}). 
Let $V$ be a real vector space of dimension $m+1$ and $n_1, \dots, n_s$ integers satisfying 
$0 < n_1 < \cdots < n_s \leq m$. Define the {\bf flag manifold} 
$\Fl(n_1, \dots, n_s, V)$ of type $(n_1, \dots, n_s)$ by 
the totality of flags $V_{n_1} \subset \cdots \subset V_{n_s} \subset V$ of linear subspaces 
with $\dim(V_i) = \ell_i, (1 \leq i \leq r)$. 
Set ${\mathcal F} = \Fl(n_1, \dots, n_s, V)$. Then the {\bf canonical distribution} $D \subset T{\mathcal F}$ is defined as follows: 
Denote by $\pi_i : {\mathcal F} \to {\Gr}(\ell_i, V)$ 
the canonical projection to the $i$-th member of the flag. 
Then, for $v \in T_{\mathbf V)}{\mathcal F}, {\mathbf V} \in {\mathcal F}$, 
$$
v \in D_{\mathbf V} \Longleftrightarrow  
{\pi_i}_*(v) \in T{\Gr}(n_i, V_{n_{i+1}}) (\subset T{\Gr}(n_i, V)), 
(1 \leq i \leq \ell-1). 
$$
Then $D$ is a subbundle of $T{\mathcal F}$ with 
$$
\rank(D) = n_1(n_2 - n_1) + (n_2 - n_1)(n_3 - n_2) + \cdots + (n_\ell - n_{\ell-1})(n - n_\ell). 
$$
Note that the flag bundle $\Fl(\ell_1, \dots, \ell_r; TP(V))$ is naturally identifies with the flag manifold 
$\Fl(1, \ell_1+1. \dots, \ell_r+1)$. Therefore the canonical differential system on $\Fl(\ell_1, \dots, \ell_r; TP(V))$ 
is induced. Then the canonical distribution $D$ on the Grassmannian bundle $\Gr(n, T(P(V)) = \Fl(1, n+1)$ 
introduced in \S \ref{Grassmannian bundle and canonical differential system} coincides 
with that introduced here. 

\bef
{\rm
Let $V$ be a real vector space of dimension $m+1$. 
A curve-germ $f : (N, a) \to P(V), \dim(N) = 1$ is called a {\bf flagged curve} 
if there exists a $D$-integral lift $\widetilde{f} : (N, a) \to \Fl(1, 2, \dots, m)$ of $f$ 
with respect to the projection $\pi_1 : \Fl(1, 2, \dots, m, V) = \Gr(1, V) = P(V)$. 
}
\enf

Let $\gamma : N \to \R P^m$ be a curve and $t_0 \in N$. 
Take a system of projective local coordinates $(x_1, x_2, \dots, x_m)$ of $\R P^m$ 
centred at $\gamma(t_0)$ and the local affine representation $(\R, t_0) \to (\R^m, 0)$, 
$$
\gamma(t) = {}^T(x_1(t), x_2(t), \dots, x_m(t))
$$ of $\gamma$. 
Consider the $(m\times k)$-matrix
$$
W_k(t) := \left( \gamma'(t_0), \gamma''(t_0), \cdots, \gamma^{(k)}(t_0)\right) 
$$
for any integer $k \geq 1$ and $k = \infty$. 
Note that the rank of $W_k(t_0)$ is independent of the choice on representations for $\gamma$. 

\bef
\label{type-def}
{\rm 
We call $\gamma$ {\bf of finite type} at $t = t_0 \in N$ if 
the $(m\times\infty)$-matrix 
$$
W_{\infty}(t_0) = \left( \gamma'(t_0), \ \gamma''(t_0), \ \cdots, \ \gamma^{(k)}(t_0), \ \cdots\cdots \ \right)
$$
is of rank $m$. 
Define, for $1 \leq i \leq m$, 
$
a_i := \min\left\{ k \mid \rank\, W_k(t_0) = i\right\}.  
$
Then we have a sequence of natural numbers $1 \leq a_1 < a_2 < \dots < a_m$, 
and we call $\gamma$ {\it of type} $(a_1, a_2, \dots, a_m)$ at $t = t_0 \in N$.  

If $(a_1, a_2, \dots, a_m) = (1, 2, \dots, m)$, then 
$t = t_0$ is called an {\it ordinary point} of $\gamma$. 
}
\enf

Let $f : N \to \R P^m$ be of finite type at $t_0 \in N$. 
Then the {\bf osculating flag} to $f$ at $t_0$ is defined by
$$
O_1(t_0) \subset O_2(t_0) \subset \cdots \subset O_k(t_0) \subset \cdots \subset O_m(t_0) = T_{f(t_0)}\R P^m, 
$$
where $O_r$ is the linear subspace of $T_{f(t_0)}\R P^m$ generated by 
$\gamma'(t_0), \ \gamma''(t_0), \ \cdots, \ \gamma^{(k)}(t_0)$. 
The corresponding projective subspace through $f(t_0)$ to $O_k(t_0)$ is also regarded
(\cite{Ishikawa95}). Then there exists unique integral lift $\widetilde{f} : N \to \Fl(1, 2, \dots, k, \dots, m, V)$ 
of $f$. 

The classification results on singularities which are related to flagged curves are given in \cite{Ishikawa10}\cite{Ishikawa12}\cite{Ishikawa12-2}.

\section{Legendre duality} 
\label{Legendre duality.} 

The Legendre duality is a natural geometric framework where the frontals play fundamental roles. 
In this section we review several studies of frontals in specified (semi-)Riemannian manifolds from \cite{Ishikawa12}\cite{IM2}. 

Let $\R^{n, m}$ denote the metric vector space of signature $(n, m)$, $n$ plus and $m$ minus 
(\cite{Harvey}\cite{O'Neill}. 

We write $\R^{n, 0}$ as $\R^n$ simply. 
Recall the space-models, the {\bf sphere} and the {\bf hyperbolic space}, 
$$
S^{n+1} = \{ x \in \R^{n+2} \mid x\cdot x = 1 \}, \quad 
H^{n+1} = \{ x \in \R^{1, n+1} \mid x\cdot x = - 1,  \ x_0 > 0 \}, 
$$
where $\R^{1, n+1} = \R^{n+2}_1 = \{ (x_0, x_1, \dots, x_{n+1}) \}$ is the Minkowski space of 
index $(1, n+1)$ (See for instance \cite{IPS}\cite{Harvey}). The inner product in 
$\R^{1, n+1}$ is defined by $x\cdot y = - x_0y_0 + \sum_{i=1}^{n+1} x_iy_i$. 
Moreover we identify Euclidean space $E^{n+1}$ with $\{ x \in \R^{n+2} \mid x_0 = 1\} \subset \R^{n+2}$ 
if necessary. 


Let $X$ denote $S^{n+1}, H^{n+1}$ or $E^{n+1}$. 
Set $Z = \widetilde{\Gr}(n, TX)$, the oriented Grassmannian bundle over $X$. 
Then $Z$ is a double covering of $\Gr(n, TX)$. 
The space $Z$ is identified with $T_1X$, the unit tangent bundle to $X$. In fact, 
$$
T_1S^{n+1} = \{ (x, y) \in S^{n+1}\times S^{n+1} \mid x\cdot y = 0 \}, 
\quad 
T_1H^{n+1} = \{ (x, y) \in H^{n+1}\times S^{1,n} \mid x\cdot y = 0 \}, 
$$
where
$S^{1, n} = \{ x \in \R^{1, n+1} \mid x\cdot x = 1\}$ is the {\bf de Sitter space}. 
Note that $Z = T_1H^{n+1}$ is identified with $T_{-1}S^{1, n} = \{ (y, v) \mid y \in S^{1, n}, v \in T_yS^{1,n}, v\cdot v = -1 \}$. 
Moreover $T_1E^{n+1} = E^{n+1}\times S^n$. 
We set $Y = S^{n+1}, S^{1, n}, \R\times S^n$ corresponding to $S^{n+1}, H^{n+1}, E^{n+1}$ respectively. 
Define $\pi_1 : Z \to X$ by the projection to the first component in three cases. Define 
$\pi_2 : Z \to Y$ by the projection to the second component in the cases $(X, Y) = (S^{n+1}, S^{n+1}), (X, Y) = (H^{n+1}, S^{1, n})$. 
In the case $(X, Y) = (E^{n+1}, \R\times S^n)$, we define $\pi_2 : Z = E^{n+1}\times S^n \to \R \times S^n$ by 
$\pi_2(x, y) = (- x\cdot y, \, y)$. 
The space has the canonical contact structure and all fibres of $\pi_1$ and $\pi_2$ are Legendre submanifold. 
Therefore $\pi_1$ and $\pi_2$ are Legendre fibrations. 
Then we have the double Legendre fibration in each case: 
$$
X \stackrel{\pi_1}{\longleftarrow} Z \stackrel{\pi_2}{\longrightarrow} Y. 
$$

As the model of duality, we do have the projective duality (\cite{Shcherbak1}\cite{IMo}): We set 
$$
Z \ = \ {\mathcal I}_{n+2} \ := \ \{ ([x], [y]) \in P^{n+1}\times P^{n+1*} \mid x\cdot y = 0 \}. 
$$
Here $P^{n+1*}$ is the dual projective space and 
$\cdot$ means the natural paring. The contact structure on ${\mathcal I}_{n+2}$ 
is defined by $dx\cdot y = x\cdot dy = 0$ (\cite{IMo}). 
The projections $\pi_1 : {\mathcal I}_{n+2} \to X = P^{n+1}, \pi_2 : {\mathcal I}_{n+2} \to Y = P^{n+1*}$ 
are both Legendre fibrations. 

The following fact is basic to unify our treatment: 

\bep
\label{flatness}
{\rm (\cite{IM}\cite{IM2})}
All Legendre double fibrations $X \longleftarrow Z \longrightarrow Y$ constructed above 
are locally isomorphic to each other. In particular each of them is locally isomorphic to 
the double fibration of the projective duality $P^{n+1} \longleftarrow {\mathcal I}_{n+2} \longrightarrow P^{n+1*}$. 
\enp
%
%


\

Let $f : N^n \to X$ be a {\it co-oriented} proper frontal (see Definition \ref{oriented-co-oriented}). 
Then there arises naturally the Legendre lift $\widetilde{f} : N \to T_1X = Z$ for $\pi_1 : Z \to X$ 
by attaching the unit normal vector field 
along $f$. The {\bf Legendre dual} of $f$ is defined by $f^\vee := \pi_2\circ \widetilde{f} : N \to Y$. 
Then $f^\vee$ is a frontal. If $f^\vee$ is a proper frontal. Then we have the equality $f^{\vee\vee} = f$. 

\

Let $\gamma : I \to X$ be a $C^\infty$ immersion from an interval or a circle $I$. 
In general, we mean by a {\bf framing} of the immersed curve 
$\gamma$, an oriented orthonormal frame 
$(e_1, e_2, \dots, e_{n+1})$ along $\gamma$.  
An immersion $\gamma$ is called {\bf framed} if a framing is given. 

\ber
{\rm 
Note that in \cite{Izumiya}\cite{Chen-Izumiya}, more general framings are considered to treat also 
light cone in Minkowski space. 
}
\enr

If $X = S^{n+1}$, then we set $e_0(t) = \gamma(t) \in S^{n+1}$, 
and we have the moving frame 
$\widetilde{\gamma} = (e_0, e_1, \dots, e_{n+1}) : I \to G = SO(n+2) \subset \GL_+(n+2, \R)$. 

If $X = H^{n+1}$, then we set $e_0(t) = \gamma(t) \in H^{n+1}$, 
and we have the moving frame 
$\widetilde{\gamma} = (e_0, e_1, \dots, e_{n+1}) : I \to G = SO(1,n+1) \subset \GL_+(n+2, \R)$. 

In any of three cases, the frame manifold $G$ is identified with an open subset of 
the oriented flag manifold $\widetilde{\mathcal F}_{n+2}$ 
consisting of oriented complete flags 
$$
V_1 \subset  V_2  \subset \cdots \subset V_{n+1} \subset \R^{n+2} 
$$
in $\R^{n+2}$. For each $g = (e_0, e_1, \dots, e_{n+1}) \in \GL_+(n+2, \R)$, 
we set the oriented subspace 
$$
V_i = \langle e_0, e_1, \dots, e_{i-1}\rangle_{\R} \subset \R^{n+2}, \ (1 \leq i \leq n+1). 
$$
This induces an open embedding $G \to \widetilde{\mathcal F}_{n+2} = \Fl(1, 2, \dots, n+1)$. 
Thus, for a framed curve $\gamma : I \to X$ in $X = E^{n+1}, S^{n+1}, H^{n+1}$, 
with the frame $(e_1, \dots, e_{n+1})$, 
we have the flagged curve $\widetilde{\gamma}$ 
by setting 
$$
V_i(t)  = \langle e_0(t), e_1(t), \dots, e_{i-1}(t) \rangle_{\R} \subset \R^{n+2}, \ (1 \leq i \leq n+1). 
$$
Then $\widetilde{\gamma}$ is a lifting of 
$\gamma$ for the projection $\pi_1 : \widetilde{\Fl}(1, 2, \dots, n+1, \R^{n+2}) \to \widetilde{\Gr}(1, \R^{n+2})$
to Grassmannian of oriented lines in $\R^{n+2}$. Note that 
there is the natural open embedding $X \subset \widetilde{\Gr}(1, \R^{n+2})$ in each of three cases. 

The projective duality plays an essential role, for instance, 
to formulate the famous Pl\" ucker-Klein's formula, 
to analyse generic projective hypersurface (Bruce, Platonova, Landis \cite{Arnold90}), 
tangent surfaces and Monge-Amp${\grave{\mbox{\rm e}}}$re equations (\cite{IM}).  

Let $f : N^{n} \lon \R{P^{n+1}}$ be a frontal. 
Then we have the Legendre lifting $\widetilde{f} : N \lon Gr(n, T\R P^{n+1}) 
= PT^*\R P^{n+1}$. 
Then we get the {\bf projective dual} $f^{\vee} : N \lon \R P^{n+1*}$ of $f$ 
by the composition of $\widetilde{f}$ 
with the projection $\pi^* : PT^*\R P^{n+1*} \lon \R P^{n+1*}$. 
If $f$ is sufficiently generic, then $f^{\vee}$ is also frontal, and we get 
the presumable equality $f^{\vee\vee} = f$. 

Viewed from Legendre duality, we consider the class of tangentially degenerate frontals. 

\bef
{\rm
Let $f : N \to \R P^{n+1} (S^{n+1}, H^{n+1}, E^{n+1})$ be a proper frontal. 
Then $f$ is called {\bf tangentially degenerate} if 
the regular locus $R(f^\vee) = \{ t \in N \mid f^\vee {\mbox{\rm \ is an immersion at\ }} t \}$ of the dual 
$f^\vee$ of $f$ is not dense in $N$. 
}
\enf

See the basic text \cite{AG} for the tangentially degenerate submanifolds.

\section{Grassmannian frontals}
\label{Grassmannian frontals}

With the notion of frontals, we are naturally led to the 
following generalization of the projective duality.
 
Let $f : N^n \lon \R P^{m}$ be a frontal of codimension 
$r = m - n$. 
Then, consider the Legendre lifting of $f$ :  
\begin{eqnarray*}
\widetilde f : N \lon \Gr(n, T\R P^{m}) 
& \hookrightarrow & \Gr(1, \R^{m+1}) \times \Gr(n+1, \R^{m+1}) \\
& \cong & \Gr(1, \R^{m+1}) \times \Gr(r, \R^{m+1*}). 
\end{eqnarray*}
The Grassmannian bundle $\Gr(n, T\R P^{m})$ is identified with 
$$
{\mathcal I} = \{ (p,q) \in \Gr(1, \R^{m+1}) \times \Gr(r, \R^{m+1*} \mid p \subseteq q^{\vee} \}. 
$$
Here, for $q \in \Gr(r, \R^{m+1*}$, we set $q^{\vee} := \{ v \in \R^{m+1} \mid \alpha(v) = 0 (\alpha \in q)\}$. 

Therefore we are naturally led to define the {\bf Grassmannian dual} 
$f^{\vee} : N \lon \Gr(r, \R^{m+1*})$ of $f : N \lon \R P^{m}$ by $\widetilde{f}$ composed with 
the projection to the second component, $(p, q) \mapsto q$. 

\bef
{\rm
A proper (co-oriented) frontal $f : N^n \to \R P^m (S^{n+1}, H^{n+1}, E^{n+1})$ 
is called {\bf tangentially degenerate} if the regular locus
$R(f^\vee) = \{ t \in N \mid f^\vee {\mbox{\rm \ is an immersion at\ }} t \}$ of the Grassmannian dual 
$f^\vee$ of $f$ is not dense in $N$. 
}
\enf

Returning to the general case, we remark that 
the equality \lq\lq\,$f^{\vee\vee} = f$\," does not have any meaning, 
even if $f^{\vee}$ is a proper frontal in the sense of Definition \ref{generalised-frontals}. 
Therefore, for a mapping into a Grassmannian, it is natural to 
specialise the definition of frontals as follows: 

Let $f : N \lon \Gr(r, \R^{m+1})$ be a $C^\infty$ mappings with $n + r \leq m + 1$. 
Set $s = m + 1 - n - r$. 
Then $f$ is called {\bf Grassmannian frontal} if 
there exists the {\it unique} integral lift $\widetilde{f} : M \lon ({\mathcal I}, {\mathcal D})$ of $f$ 
with respect to a fibration $\pi_1 : {\mathcal I} \lon \Gr(r, \R^{m+1})$ and a distribution ${\mathcal D}$ 
on ${\mathcal I}$ defined as follows:  
First set   
$$
{\mathcal I} := \{ (p,q) \in \Gr(r, \R^{m+1}) \times \Gr(s, \R^{m+1*}) \mid 
p \subseteq q^{\vee} \},  
$$
and consider the projection $\pi_1 : I \to \Gr(r, \R^{n+2})$ 
(resp.  $\pi_2 : I \to \Gr(s, \R^{m+1*})$). 
Moreover set 
$$
{\mathcal P} := \{ (p,q,p') \in \Gr(r, \R^{m+1}) \times \Gr(s, \R^{m+1*}) \times \Gr(r, \R^{m+1}) \mid 
p \subseteq q^{\vee} \ , p' \subseteq q^{\vee} \}, 
$$
and consider the projection $\rho : {\mathcal P} \to {\mathcal I}$ to the first and second factors 
(resp. $\varphi : {\mathcal P} \to \Gr(r, \R^{m+1})$ to the third factor).  
Then we get the double fibration $(\rho, \ \varphi)$:
$$
{\mathcal I} \stackrel{\rho}{\longleftarrow} {\mathcal P} \stackrel{\varphi}{\lon} \Gr(r, \R^{m+1}).  
$$
For each $c = (p, q) \in {\mathcal I}$, we consider $\rho^{-1}(c)$. 
Then we consider its projection 
$$
\varphi(\rho^{-1}(c)) = 
\{p' \in \Gr(r, \R^{m+1}) \mid p' \subseteq q^{\vee}\} 
$$ 
by $\varphi$, which is regarded as $\Gr(r, \R^{r+n})$. 
Note that $\dim q^{\vee} = r+n, \ p \in \varphi(\rho^{-1}(c))$ 
and that $\varphi(\rho^{-1}(c)) \subset \Gr(r, \R^{m+1})$ is a 
submanifold of codimension $r(m+1-r) - rn = rs$. 

Define the {\it tautological subbundle} ${\mathcal D} \subset T{\mathcal I}$ 
of codimension $rs$, for each $c = (p,q) \in {\mathcal I}$,  
by 
$$
{\mathcal D}_c = \pi_*^{-1}(T_p(\varphi(\rho^{-1}(c)))) \subset T_c{\mathcal I}. 
$$ 
Note that, if $r \not= 1$, or, $r \not = n + 1$, 
then the \lq\lq system of tangential linear subspaces\rq\rq\,  
$\{ \varphi(\rho^{-1}(c)) \mid c \in I \}$ 
in the Grassmannian $\Gr(r, \R^{m+1})$ 
defined by ${\mathcal D}$ does not represent general tangential linear subspaces 
of the Grassmannian. 

If we take local Grassmannian coordinates $(a_{ij})_{1\leq i \leq r, 1 \leq j \leq n+s}$  
of $\Gr(r, \R^{m+1})$ and 
$(b_{k\ell})_{1\leq k \leq n+r, 1\leq \ell \leq s}$ of $\Gr(s, \R^{m+1*})$, then 
${\mathcal I}$ is defined by the system of equations 
$$
b_{ij} + a_{i1}b_{r+1\, j} + \cdots + a_{in}b_{r+n\, j} + a_{i\, n+j} = 0, 
\ 1 \leq i \leq r, 1 \leq j \leq s,
$$
and ${\mathcal D}$ is defined by the system of $1$-forms 
$$
b_{r+1\, j}da_{i1} + \cdots + b_{r+n\, j}da_{in} + da_{i\, n+j} = 0, 
\ 1 \leq i \leq r, 1 \leq j \leq s. 
$$

The integral lifting $\widetilde{f}$ 
is called the {\bf Legendre lifting} of $f$ in the generalised sense. 
The relation to the original definition of frontals is as follows: 

\bel
Let $F : (\R^n, 0) \lon ({\mathcal I},(p_0,q_0))$ 
be an integral map-germ to the distribution ${\mathcal D} \subset T{\mathcal I}$. 
Then $f = \pi_1\circ F : (\R^n, 0) \lon (\Gr(r, \R^{m+1}), p_0)$ is 
Grassmannian frontal if and only if $\kappa\circ f$ is proper, i.e. 
$S(\kappa\circ f) \subset (\R^n, 0)$ is nowhere dense, for 
some projection 
$$
\kappa : (\Gr(r, \R^{m+1}),p_0) \hookrightarrow 
(Hom(\R^r, \R^{n+s}),0) \stackrel{i^*}{\lon} (Hom(\R, \R^{n+s}),0)  
\hookrightarrow \R P^{n+s-1}, 
$$
induced from a linear inclusion $i : \R \hookrightarrow \R^r$. 
\enl 

Now, from the duality, we have another 
distribution ${\mathcal D}' \subset T{\mathcal I}$ 
from the projection $\pi' : {\mathcal I} \lon \Gr(s, \R^{m+1*})$ to the second factor, 
setting 
$$
{\mathcal P}' = \{ (q', p, q) \in \Gr(s, \R^{m+1*}) \times 
\Gr(r, \R^{n+2}) \times \Gr(s, \R^{m+1*}) \mid 
q \subseteq p^{\vee}, \ q' \subseteq p^{\vee} \}. 
$$
Then the fundamental result is the following: 

\bep
Two distributions ${\mathcal D}$ and ${\mathcal D}'$ on the incidental manifold ${\mathcal I}$ coincide. 
\enp

We conclude this section by the following observation: 

\bep
Let $F : N^n \to I \subset Gr(r, \R^{m+1}) \times Gr(s, \R^{m+1*})$ be an integral mapping to 
the distribution ${\mathcal D}$ with $n + r + s = m+1$. 
Suppose $\pi\circ F =: f$ and $\pi'\circ F =: f\vee$ are Grassmannian frontals respectively. 
Then we have $f^{\vee\vee} = f$. 
\enp

\section{Tangent varieties}
\label{Tangent varieties}

Given a curve in Euclidean $3$-space $\EE^3 = \R^3$, 
the embedded tangent lines to the curve draw a surface in $\R^3$, 
which is called the {\bf tangent surface} (or {\bf tangent developable}) to the curve (\cite{Cayley}\cite{Ishikawa96}). 
It is known that the tangent surfaces (tangent developables) are developable surfaces. 
Developable surfaces which are locally isometric to the plane keep on interesting 
many mathematicians, for instance, 
Monge (1764), Euler (1772), Cayley (1845), Lebesgue (1899). See \cite{Lawrence} for details. 
Therefore the tangent surfaces are regarded as generalised solutions (with singularities) 
of the Monge-Amp\`{e}re equation 
$$
\dfrac{\pa^2 z}{\pa x^2}\dfrac{\pa^2 z}{\pa y^2} - \left(\dfrac{\pa^2 z}{\pa x\pa y}\right)^2 = 0
$$
on spacial surfaces $z = z(x, y)$. 
Tangent surfaces are flat in $E^3$. However they are not flat but \lq\lq extrinsically flat" or tangentially degenerate 
in $S^3$ and $H^3$ (cf. \cite{AG}\cite{KRSUY}). See also \S \ref{Legendre duality.}. 
The notion of {\bf types} $(a_1, a_2, a_3)$ for a curve-germ is introduced (Definition \ref{type-def}). Then 
the cuspidal edge, (resp. the swallowtail, the cuspidal beaks (Mond surface), the cuspidal butterfly) is obtained 
as the tangent developable of a curve of type $(1, 2, 3)$ (resp. $(2, 3, 4)$, $(1, 3, 4)$, $(3, 4, 5)$). 

This property is related to \lq\lq projective duality": 
{\it The projective dual of a tangent surface collapse to a curve  (the dual curve)}. See \cite{Ishikawa99}. 

Let $\gamma : \R \to \R^3$ be an immersed curve. 
Then the tangent surface has the natural parametrization 
$$
\Tan(\gamma) : \R^2 \to \R^3, \quad \Tan(\gamma)(t, s) := \gamma(t) + s\gamma'(t). 
$$
The tangent surface necessarily has singularities at least along $\gamma$, 
\lq\lq the edge of regression". 

It is known that the tangent surface to a generic curve $\gamma : \R \to \R^3$ in $\R^3$ 
has singularities only along $\gamma$ and is locally diffeomorphic to the cuspidal edge 
or to the folded umbrella (also called, the cuspidal cross cap), 
as is found by Cayley and Cleave (1980). 
Cuspidal edge singularities appear along ordinary points where 
$\gamma', \gamma'', \gamma'''$ are linearly independent, while 
the folded umbrellas appear at isolated points of zero torsion where 
$\gamma', \gamma'', \gamma'''$ are linearly dependent but $\gamma', \gamma'', \gamma''''$ 
are linearly independent. 

In a higher dimensional space $\R^m, m \geq 4$, 
for an immersed curve $\gamma : \R \to \R^m$, we define 
the tangent surface 
$\Tan(\gamma) : \R^2 \to \R^m$ \ 
by $\Tan(\gamma)(t, s) := \gamma(t) + s\gamma'(t)$. 
Then we have generically that 
$\gamma', \gamma'', \gamma'''$ are linearly independent and 
$\Tan(\gamma)$ is locally diffeomorphic to 
the (embedded) cuspidal edge in $\R^m$. 
Now we give the general definition: 

\

\bef
{\rm 
Let $N$ be an $n$-dimensional manifold. Let $f : N^n \to \R^m$ be a proper frontal. 
Let $\widetilde{f} : N \to \Gr(n, T\R^m)$ be the Legendre lift of $f$. 
Then the {\bf tangent mapping} 
$\Tan(f) : T_f \to \R^m$ of $f$ is defined by, for $t \in N$ and $v \in \widetilde{f}(t) \subset T_{f(t)}\R^m$, 
$$
\Tan(f)(t, v) := f(t) + v, \quad (t, v) \in T_f, 
$$
using the affine structure of $\R^m$. 
Then we define the {\bf tangent variety} of $f$ as the parametrised variety which is 
defined by the right equivalence class of $\Tan(f)$. 
If $(t_1, \dots, t_n)$ is a system of local coordinates of $N$, and 
$(t_1, \dots, t_n, s_1, \dots, s_n)$ the induced system of local coordinates of 
$T_f$ induced by a system of local frame $v_1(t), \dots, v_n(t)$ of 
$\widetilde{f}$, then $\Tan(f)$ is given by 
$$
\Tan(f)(t, s) = f(t) + \sum_{j=1}^n s_j v_j(t)). 
$$
}
\enf

Also note that we can define similarly the tangent varieties of mappings to a projective space. 
Tangent varieties appear in various geometric problems and applications naturally (\cite{AG}\cite{BG}\cite{IMT1}\cite{IMT2}\cite{IMT3}\cite{IKY}\cite{IS}\cite{INS}\cite{Porteous}\cite{Lawrence}\cite{Zak}). 
See \cite{Ishikawa99}\cite{Ishikawa12-2}, for the geometric exposition on the local classification problem of tangent varieties. In particular it is proved in \cite{Ishikawa12-2}\cite{IMT3} 
the following: 

\bep
Let $\gamma : (N, t_0) \to \R P^m$ be a curve-germ of finite type (Definition \ref{type-def}). 
Then $\Tan(\gamma) : (N\times\R, (t_0, 0)) \to \R P^m$ is a proper frontal. 
\enp

A proper frontal $f : N \to M$ is called a {\bf directed curve} if $\dim(N) = 1$ (\cite{IY}\cite{IY1}\cite{IY2}). 
A directed curve $\gamma$ is called {\bf orientable} if 
there exists a frame $u : N \to TM$, $u(t) \not= 0$, along $\gamma$ such that 
$
\gamma'(t) \in \langle \ u(t) \ \rangle_{\R}, \ t \in \R, 
$
which projects to the unique Legendre lift $\widetilde{\gamma} : N \to P(TM) = \Gr(1, TM)$ of $\gamma$ 
satisfying 
$
\gamma'(t) \ \in \ \widetilde{\gamma}(t), \ (t \in \R). 
$

Let $\gamma : N \to M$ be a directed curve and $\widetilde{\gamma} : N \to P(TM)$ the unique $D$-integral lift of $f$. 
Recall that the {\bf tangent bundle} to $f$ is defined by 
$
T_\gamma := \{(t, v) \in N \times TM \mid v \in \widetilde{\gamma}(t)\}, 
$
which is a line bundle over $N$ (see \S \ref{tangent-and-normal}). 
Let $M$ be a manifold with an affine connection. 
We define the {\bf tangent mapping} $\Tan(\gamma) : T_\gamma \to M$ by $(t, v) \to \exp(v)$, 
using the exponential map (see \S \ref{Affine connection and tangent surface}). 

\ber
{\rm
By Lemma \ref{global-n=1}, there exists a global Legendre lift $\widetilde{\gamma} : N \to P(TM)$ of $f$. 
Then the orientability condition means that the line bundle 
$T_{\widetilde{f}}$ over $N$ is orientable. 
}
\enr

Let $M = \R^m$ and $\gamma : N \to \R^m$ be a directed and orientable curve. 
Then the tangent surface $\Tan(\gamma) : N\times \R \to \R^m$ of a directed curve $\gamma$ is defined by 
$$
\Tan(\gamma)(t, s) := \gamma(t) + s \, u(t)
$$ 
The right equivalence class of $\Tan(\gamma)$ is independent of the choice of frame $u$. 

The singularities of the tangent surface $\Tan(\gamma)$ 
for a generic directed curve $\gamma : \R \to \R^m$ 
on a neighbourhood of the curve are only the cuspidal edge, the folded umbrella, 
and swallowtail if $m = 3$, and 
the embedded cuspidal edge and the open swallowtail if $m \geq 4$. See \cite{Cleave}\cite{IY}. 
Several degenerate cases are studied in 
\cite{Mond1}\cite{Mond2}\cite{Ishikawa93}\cite{Ishikawa95}\cite{Ishikawa96}\cite{Ishikawa99}.

\

\section{Grassmannian geometry}
\label{Grassmannian geometry}

We will give a series of classification results of singularities of tangent surfaces 
in $A_n$-geometry, i.e. the geometry associated to the group 
${\mbox{\rm PGL}}(n+1, \R)$ (see \cite{IMT4}). 

Let $V = \R^{m+1}$ be the vector space of dimension $m+1$ and consider 
a flag in $V$ of the following type (a complete flag): 
$$
V_1 \subset V_2 \subset V_3 \subset \cdots \subset V_m \subset V, \quad \dim(V_i) = i. 
$$
The set of such flags form a manifold ${\Fl}(1,2,3,\dots, m)$ 
of dimension $\frac{n(n+1)}{2}$. 

A curve $\gamma : \R \to P(V) = P(V^{m+1})$ 
arises a $D$-integral curve $\Gamma : \R \to {\Fl}(1,2,3,\dots, m)$ for the canonical 
distribution $D \subset T{\Fl}(1,2,3,\dots, m)$, if we regard its osculating planes: 
the curve itself is given by $V_1(t)$, the tangent line is given by 
$V_2(t)$, the osculating plane is given by $V_3(t)$ and so on. 

Let $m = 2$. Let $V_1(t) \subset V_2(t) \subset V = \R^3$ be an admissible curve. 
For each $a$, planes $V_2$ satisfying $V_1(a) \subset V_2 \subset V$ 
form the tangent line to the curve $\{ V_1(t)\}$ at $t = a$ in $P(V) = 
P^2$. Similarly lines $V_1$ satisfying $V_1 \subset V_2(a)$ form the 
tangent line to the dual curve $\{ V_2(t)\}$ at $t = a$ in $\Gr(2, V) = 
P(V^*) = P^{2*}$, the dual projective plane. 
For a generic admissible curve, we have the duality on 
\lq\lq tangent maps\rq\rq: 

Let $m = 3$. 
Let $\Gamma : \R \to \Fl(1, 2, 3)$ be a $D$-integral curve. Set 
$\Gamma(t) = (V_1(t), V_2(t), V_3(t))$, 
$V_1(t) \subset V_2(t) \subset V_3(t) \subset V = \R^4$. Then $\Gamma$ 
induces the curve $\pi_1\circ \Gamma$ in $P^3 = P(\R^4)$, the curve $\pi_2\circ\Gamma$ in $\Gr(2, \R^4)$ 
and the curve $\pi_3\circ\Gamma$ in $P^{3*} = \Gr(3, \R^4)$. 
Then we have the following duality on their tangent surfaces in $A_3$-geometry: 
$$
\begin{array}{c|c|c}
{\mbox{\rm \textcircled{\scriptsize 3}}} & {\mbox{\rm \textcircled{\scriptsize 4}}} & {\mbox{\rm \textcircled{\scriptsize 3}}}
\\
\hline
{\mbox{\rm 
Cuspidal \ Edge}} & {\mbox{\rm Cuspidal \ Edge}} & {\mbox{\rm Cuspidal \ Edge}} \\
\hline
{\mbox{Swallow Tail}} & {\mbox{\rm Cuspidal \ Edge}} & {\mbox{Folded \ Umbrella}}\\
\hline
{\mbox{Mond \ Surface}} & {\mbox{Open \ Swallowtail}} & {\mbox{ Mond \ Surface}}
\\
\hline
{\mbox{Folded \ Umbrella}} & {\mbox{Cuspidal \ Edge}} & {\mbox{Swallow Tail}}
\end{array}
$$\ 

In general for a generic $D$-integral curve $\Gamma : \R \to \Fl(1, 2, 3, \dots, m)$, 
$$
V_1(t) \subset V_2(t) \subset V_3(t) \subset \cdots \subset V_m(t) \subset V = \R^{m+1}, 
$$
we have the classification of singularities of tangent surfaces:

\bet 
{\rm ($A_n, n \geq 4$)} 
The classification list consists of $n+1$ cases for curves in Grassmannians:
\begin{center}
\begin{tabular}{|cccccc|}
\hline
$P^n$ & $\Gr(2, V)$ & $\Gr(3, V)$ & $\Gr(4, V)$  &$\cdots$ & $\Gr(n, V)$
\\
\hline
\hline
{\rm CE} & {\rm CE} & {\rm CE} & {\rm CE} & $\cdots$ & {\rm CE} 
\\
\hline
{\rm OSW} & {\rm CE} & {\rm CE} & {\rm CE} & $\cdots$ & {\rm CE} 
\\
\hline
{\rm OM} & {\rm OSW} & {\rm CE} & {\rm CE} & $\cdots$ & {\rm OFU} 
\\
\hline
{\rm OFU}  & {\rm CE} & {\rm OSW} & {\rm CE} & $\cdots$ & {\rm OM} 
\\
\hline
{\rm CE} & {\rm CE} & {\rm CE} & {\rm OSW} & $\cdots$ & {\rm CE} 
\\
\hline
$\vdots$ & $\vdots$ & $\vdots$ & $\vdots$ & {$\ddots$} & $\vdots$ 
\\
\hline
{\rm CE} & {\rm CE} & {\rm CE} & {\rm CE} & $\cdots$ & {\rm OSW}
\\
\hline
\end{tabular}
\end{center}
\ent

The {\bf cuspidal edge} (resp. {\bf open swallowtail}, {\bf open Mond surface}, 
{\bf open folded umbrella})
is defined as a diffeomorphism class of the tangent surface-germ to a curve of type
$(1,2,3,\cdots)$ (resp. $(2,3,4,5,\cdots)$, $(1,3,4,5,\cdots)$, $(1,2,4,5,\cdots)$) in an affine space. 

\


\section{Affine connection and tangent surface} 
\label{Affine connection and tangent surface} 

Now let us consider the case of directed curves in a Riemannian manifold, 
or more generally, the case of directed curves in 
a manifold with any affine connection, which is not necessarily projectively flat. 
For any directed curve, we have the well-defined tangent geodesic to each point of the curve. 
If we regard it as the \lq\lq tangent line", then we have the well-defined 
tangent surface for the directed curve.

It is proved in \cite{IY}, 
for any affine connection on a manifold of dimension $m \geq 3$, 
the singularities of the tangent surface to a generic directed curve 
on a neighbourhood of the curve are only the {\bf cuspidal edge}, the {\bf folded umbrella}, 
and {\bf swallowtail} if $m = 3$, and 
the {\bf embedded cuspidal edge} and the {\bf open swallowtail} if $m \geq 4$. 
Moreover we have: 

\bet
\label{Characterisation-tangent}
{\rm (\cite{IY})} 
Let $\nabla$ be any torsion-free affine connection on a manifold $M$. 
Let $\gamma : \R \to M$ be a $C^\infty$ curve. 

{\rm (1)}  Let $\dim(M) = 3$. 
If $(\nabla\gamma)(a), (\nabla^2\gamma)(a), (\nabla^3\gamma)(a)$ are linearly independent at $t = a \in \R$, 
then the tangent surface $\Tan(\gamma)$ is locally diffeomorphic to the cuspidal edge at $(a, 0) \in \R^2$. 
If $(\nabla\gamma)(a), (\nabla^2\gamma)(a), (\nabla^3\gamma)(a)$ are linearly dependent, 
and $(\nabla\gamma)(a), (\nabla^2\gamma)(a), (\nabla^4\gamma)(a)$ are linearly independent, then
the tangent surface $\Tan(\gamma)$ is locally diffeomorphic to the folded umbrella at $(a, 0) \in \R^2$. 
If $(\nabla\gamma)(a) = 0$ and $(\nabla^2\gamma)(a), (\nabla^3\gamma)(a), (\nabla^4\gamma)(a)$ 
are linearly independent, then
the tangent surface $\Tan(\gamma)$ is locally diffeomorphic to the swallowtail at $(a, 0) \in \R^2$. 
 
{\rm (2)} Let $\dim(M) \geq 4$. 
If $(\nabla\gamma)(a), (\nabla^2\gamma)(a), (\nabla^3\gamma)(a)$ are 
linearly independent at $t = a \in \R$, 
then the tangent surface $\Tan(\gamma)$ is locally diffeomorphic to the 
embedded cuspidal edge at $(a, 0) \in \R^2$. 
If $(\nabla\gamma)(a) = 0$ and 
$
(\nabla^2\gamma)(a), (\nabla^3\gamma)(a), (\nabla^4\gamma)(a), (\nabla^5\gamma)(a)
$ 
are linearly independent at $t = a \in \R$, 
then the tangent surface $\Tan(\gamma)$ is locally diffeomorphic to the 
open swallowtail at $(a, 0) \in \R^2$. 
\ent

For the proof of Theorem \ref{Characterisation-tangent}, we apply the characterisation theorems found in \cite{KRSUY}\cite{FSUY}\cite{Ishikawa12-2}. 

In \cite{IY}\cite{IY1}, singularities of tangent surfaces of torsionless curves are studied. In that case, so called fold singularities and $(2, 5)$-cuspidal edges appear. See also \cite{HKS}. 
%
%

\section{Characterisation of frontal singularities}
\label{Characterisation of frontal singularities}

When we treat singularities in a general ambient space as in the previous section, we need the intrinsic 
characterisations of singularities. 
Note that the characterization of swallowtails was applied to hyperbolic geometry in \cite{KRSUY} 
and to Euclidean and affine geometries in \cite{IM}. 
The characterization of folded umbrellas is applied to Lorenz-Minkowski geometry in \cite{FSUY}. 
In Theorem \ref{Characterisation-tangent}, we apply to non-flat projective geometry the characterisations and their some generalization 
via the notion of openings introduced in \S \ref{Openings and frontals}. 

Let $f : (\R^2, p) \to M^3$ be a frontal with a non-degenerate singular point at $p$ (see Lemma \ref{non-degenerate-proper}) 
and $\widetilde{f} : (\R^2, p) \to \Gr(2, TM)$ the integral lifting of $f$. 
Let $V_1, V_2 : (\R^2, p) \to TM$ be an associated frame with $\widetilde{f}$. 
Let $L : (\R^2, p) \to T^*M \setminus \zeta$ be an annihilator of 
$\widetilde{f}$. The condition is that $\langle L, V_1\rangle = 0, \langle L, V_2\rangle = 0$. 
Here $\zeta$ means the zero section. 
Let $c : (\R, t_0) \to (\R^2, p)$ be a parametrization of the singular locus $S(f)$, $p = c(t_0)$,  and 
$\eta : (\R^2, p) \to T\R^2$ be a vector field which restricts to the kernel field of $f_*$ on $S(f)$. 
Suppose that $V_2(p) \not\in f_*(T_p\R^2)$. 
Then, for any affine connection $\nabla$ on $M$, we define 
$$
\psi(t) := \langle L(c(t)),  (\nabla_\eta^f V_2)(c(t))\rangle. 
$$
Note that the vector field $(\nabla_\eta^f V_2)(c(t))$ is independent of the extension $\eta$ 
and the choice of affine connection $\nabla$, since $\eta\vert_{S(f)}$ is a kernel field of $f_*$. 
We call the function $\psi(t)$
the {\bf characteristic function} of $f$.

Then the following characterisations of cuspidal edges and folded umbrellas are given in \cite{KRSUY}\cite{FSUY}: 

\bet
\label{characterization}
{\rm (Theorem 1.4 of \cite{FSUY})}. 
Let $f : (\R^2, p) \to M^3$ be a germ of frontal with a non-degenerate singular point at $p$. 
Let $c : (\R, t_0) \to (\R^2, p)$ be a parametrization of the singular locus of $f$. 
Suppose $f_*c'(t_0) \not= 0$. Then, for the characteristic function $\psi$, 
\\
{\rm (1)} 
$f$ is diffeomorphic to the cuspidal edge if and only if $\psi(t_0) \not= 0$. 
\\
{\rm (2)} $f$ is diffeomorphic to the folded umbrella if and only if
$\psi(t_0) = 0, \psi'(t_0) \not= 0$. 
\ent

We can summarise several known results as those on openings of the fold: 

\bet
\label{characterization2}
{\rm (\cite{IY}\cite{IY1})}
Let $f : (\R^2, p) \to M^m, m \geq 2$ be a germ of frontal with a non-degenerate singular point at $p$, 
$\widetilde{f} : (\R^2, p) \to \Gr(2, TM)$ the integral lifting of $f$ and 
$V_1, V_2 : (\R^2, p) \to TM$ an associated frame with $\widetilde{f}$. 
Let $c : (\R, t_0) \to (\R^2, p)$ be a parametrization of the singular locus of $f$. 
Suppose $f_*c'(t_0) \not= 0$. 
Then $f$ is diffeomorphic to an opening of the fold, namely to the germ $(u, w) \mapsto (u, \frac{1}{2}w^2)$. 
Moreover we have: 
\\
{\rm (0)} Let $m = 2$. Then $f$ is diffeomorphic to the fold. 
\\
{\rm (1)} Let $m \geq 3$. Then $f$ is diffeomorphic to the cuspidal edge if and only if
$\psi(t_0) \not= 0$. 
\\
{\rm (2)} Let $m = 3$. Then $f$ is diffeomorphic to the folded umbrella if and only if
$\psi(t_0) = 0, \psi'(t_0) \not= 0$. 
\ent

Based on results in \cite{KRSUY} and \cite{Ishikawa12-2}, 
we summarise the characterization results on openings of the Whitney's cusp map-germ: 

\bet
\label{characterization3}
{\rm (\cite{IY}\cite{IY2})}
Let $f : (\R^2, p) \to M^m, m \geq 2$ be a germ of frontal with a non-degenerate singular point at $p$, 
$V_1, V_2 : (\R^2, p) \to TM$ an associated frame with $\widetilde{f}$ with $V_2(p) \not\in f_*(T_p\R^2)$, 
and $\eta : (\R^2, p) \to T\R^2$ an extension of a kernel field along of $f_*$. 
Let $c : (\R, t_0) \to (\R^2, p)$ be a parametrization of the singular locus of $f$. 
Set $\gamma = f\circ c : (\R, t_0) \to M$. 
Suppose $(\nabla\gamma)(t_0) = 0$ and $(\nabla^2\gamma)(t_0) \not= 0$. 
Then $f$ 
is diffeomorphic to an opening of Whitney's cusp, 
namely to the germ $(u, t) \mapsto (u, t^3 + ut)$. Moreover we have 

{\rm (0)} Let $m = 2$. Then $f$ is diffeomorphic to Whitney's cusp. 

{\rm (1)} Let $m = 3$. Then $f$ is diffeomorphic to the swallowtail if and only if
$$
V_1(c(t_0)), \ V_2(c(t_0)), \ (\nabla^f_\eta V_2)(c(t_0))
$$ 
are linearly independent in $T_{f(p)}M$. 

{\rm (2)} Let $m \geq 4$. Then $f$ is diffeomorphic to the open swallowtail if and only if
$$
(V_1\circ c)(t_0), \ (V_2\circ c)(t_0), \ ((\nabla^f_\eta V_2)\circ c)(t_0), \ 
(\nabla^\gamma_{\pa/\pa t}((\nabla^f_\eta V_2)\circ c))(t_0)
$$ 
are linearly independent in $T_{f(p)}M$.  
\ent

Note that the conditions appeared in Theorem \ref{characterization} are invariant under diffeomorphism equivalence 
introduced in Introduction. 
In fact the conditions are invariant under a weaker equivalence relation. 
In Definition \ref{J-equivalence}, we have introduce the notion of ${\mathcal J}$-equivalence 
of map-germs. 

\bec
Let $f : (\R^2, 0) \to (\R^m, 0)$ be a frontal. Then $f$ is ${\mathcal J}$-equivalent to Whitney's cusp 
if and only if $f$ is diffeomorphic to an opening of Whitney's cusp. Moreover. if $m = 2$, then 
$f$ is diffeomorphic to Whitney's cusp. If $m = 3$ and $f$ is a front, then $f$ is diffeomorphic to swallowtail. 
\enc

The known criteria of singularities (see for instance \cite{Saji10}\cite{Saji11}\cite{Kabata}) seem to be closely related 
with frontals, openings and ${\mathcal J}$-equivalence. The detailed relations are still open to be studied.

\section{Null frontals}
\label{Null frontals}

Let $(M, g)$ be a semi-Riemannian manifold with an indefinite metric $g$. 
Denote by ${\mathcal C} \subset TM$ the null cone field associated with the indefinite metric $g$, i.e. 
${\mathcal C}$ is the set of null vectors: 
$$
{\mathcal C} = \bigcup_{x \in M} {\mathcal C}_x, \quad {\mathcal C}_x = \{ u \in T_xM \mid \ g_x(u, u) = 0 \}. 
$$
Let $\pi : {\mathcal C} \to M$ be the canonical projection. 

\bef
{\rm 
A mapping $f : N \to M$ is called {\bf totally null} (resp. {\bf null}) 
if the induced metric $f^*g$ is identically zero (resp. $f^*g$ is degenerate everywhere). 
The condition that $f$ is totally null is equivalent to that $f_*(T_tN) \subset {\mathcal C}_{f(t)}$ 
(resp. $f_*(T_tN)$ is tangent to ${\mathcal C}_{f(t)}$), for any $t \in N$. 
}
\enf

\bef
{\rm 
A curve-germ $\gamma : (\R, a) \to M$ is called {\bf null} if 
$
\gamma'(t) \in {\mathcal C} \quad  (t \in (\R, a)). 
$
Moreover $\gamma : (\R, a) \to M$ is called {\bf null-directed} if there exists a lift  
$u : (\R, a) \to {\mathcal C}$ such that 
$\pi\circ u = \gamma, u(t) \not= 0, \gamma'(t) \in \langle u(t) \rangle_\R, \ t \in (\R, a). $
}
\enf

A map-germ is null (resp. null-directed) if and only if it is totally null (resp. totally null frontal). 

\bef
{\rm 
Let $\gamma : (\R, a) \to M$ be null-directed. 
Define the {\bf null tangent surface} $\Tan(\gamma) : (\R^2, a\times \R) \to M$ of $\gamma$ as the 
ruled surface by null geodesics through points $\gamma(t)$ with the directions $u(t)$. 
}
\enf

The right equivalence class of $\Tan(\gamma)$ is independent of the choice of the lift $u$. 

We have the following classification results. For the details see \cite{IMT3}\cite{IMT4}: 
The singularities of tangent surface $\Tan(\gamma)$ for a generic null directed curve 
$\gamma : \R \to \R^{2,2}$ are cuspidal edges and open swallowtails. 
The singularities of tangent surface $\Tan(\gamma)$ for a generic null directed curve 
$\gamma : \R \to \R^{2,3}$ are  cuspidal edges, open swallowtails, 
open Mond surfaces and unfurled folded umbrellas. 
The singularities of tangent surface $\Tan(\gamma)$ for a generic null directed curve 
$\gamma : \R \to \R^{3,3}$ (the projection of a generic \lq\lq Engel integral\rq\rq\ curve) 
are embedded cuspidal edges, open swallowtails and open Mond surfaces. 
See \cite{Ishikawa12-2} for the normal forms and pictures of the singularities. 

\

In general the tangent surface 
to a null curve is a ruled surface by null lines, which is not necessarily 
a totally null surface, but a null surface, which we call the {\it null tangent surface}.

Let $X$ be a $3$-dimensional Lorentzian manifold {\rm(}with signature $(1, 2)${\rm )}.  
A smooth map-germ $\varphi : (\R^2, 0) \to X$ is called a {\bf null frontal surface} 
or a {\bf null frontal} in short if 
there exists a smooth lift $\widetilde{\varphi} : (\R^2, 0) \to PT^*X = \Gr(2, TX)$ of $\varphi$ 
such that 
$\widetilde{\varphi}(t)$ is a lightlike plane in $T_{\varphi(t)}X$ 
and $\varphi_*(T_t\R^2) \subset \widetilde{\varphi}(t)$, for any $t \in (\R^2, 0)$. 
The notion of null frontals is a natural generalization of null immersions to singular surfaces. 
We have presented several classification results of singularities which arise in null frontals up to local diffeomorphisms and up to $O(2, 3)$-conformal transformations in the conformally flat case {\rm (}cf. \cite{IMT1}{\rm )}. 
The classification is achieved by using the fact that null frontals are obtained as tangent surfaces to null curves in $X$, as well as 
\lq\lq associated varieties\rq\rq\, to Legendre curves in the space $Y$ 
of null geodesics on $X$ {\rm (}cf. \cite{IMT2}\cite{IMT3}\cite{IMT4}{\rm )}. 
A related result is obtained in \cite{Chino-Izumiya}. 


%
%
\section{Abnormal frontals}
\label{Abnormal frontals}

Let $M$ be a $5$-dimensional manifold and 
${\mathcal D} \subset TM$ a distribution of rank $2$. 
Then ${\mathcal D}$ is called a {\bf Cartan distribution} if it has growth $(2, 3, 5)$, namely, if 
$\rank({\mathcal D}^{(2)}) = 3$ and $\rank({\mathcal D}^{(3)}) = 5$, where, 
we define in terms of Lie bracket, ${\mathcal D}^{(2)} = 
{\mathcal D} + [{\mathcal D}, {\mathcal D}]$ and ${\mathcal D}^{(3)} = 
{\mathcal D}^2 + [{\mathcal D}, {\mathcal D}^2]$. 
It is known that, for any point $x$ of $M$ and for any direction $\ell \subset {\mathcal D}_x$,  
there exists an abnormal geodesic, which is unique up to parametrisations, 
through $x$ with the given direction $\ell$ (see \cite{IKY1}\cite{IKY2}). 

Then, for a given ${\mathcal D}$-directed curve $\gamma$, 
we define {\bf abnormal tangent surface} of $\gamma$, which is ruled by
abnormal geodesics through points $\gamma(t)$ with the directions $u(t)$. 

On $\R^5$ with coordinates $(\lambda, \nu, \mu, \tau, \sigma)$, 
define the distribution ${\mathcal D} \subset T\R^6$ generated by the pair of vector fields 
$$
\begin{array}{rcl}
\eta_1 & = & \dfrac{\pa}{\pa \lambda} + \nu\dfrac{\pa}{\pa \mu} 
- (\lambda\nu - \mu)\dfrac{\pa}{\pa \tau} + \nu^2\dfrac{\pa}{\pa \sigma}, 
\vspace{0.3truecm}
\\
\eta_2 & = & \dfrac{\pa}{\pa \nu} - \lambda\dfrac{\pa}{\pa \mu} 
+ \lambda^2\dfrac{\pa}{\pa \tau} 
- (\lambda\nu + \mu)\dfrac{\pa}{\pa \sigma}. 
\end{array}
$$
Then ${\mathcal D} \subset T\R^6$ is a Cartan distribution and it has maximal symmetry of dimension $14$, 
maximal among all Cartan distributions, 
which is of type $G_2$, one of simple Lie algebras. 

For a generic $G_2$-Cartan directed curve $\gamma : \R \to 
\R^5$, the tangent surfaces at any point $a \in \R$ 
is classified, up to local diffeomorphisms, into 
{\bf embedded cuspidal edge}, {\bf open Mond surface}, and {\bf generic open folded pleat} (see \cite{IMT2} for details). 
The classification of singularities in abnormal tangent surfaces to generic Cartan directed curves 
for general Cartan distributions seems to be un-known yet. 

\

\section{Appendix: Malgrange preparation theorem on differentiable algebras}
\label{Appendix: Malgrange preparation theorem on differentiable algebras}

We show the Malgrange's preparation theorem on differentiable algebras \cite{Malgrange}
from the ordinary Malgrange-Mather's preparation theorem (see for example \cite{Brocker}), relating to the theory of 
$C^\infty$-rings which we have utilised in this paper. 

An $\R$-algebra $A$ is called {\bf local} if 
it has a unique maximal ideal ${\mathfrak m}_A$.  

\bee
\label{E_n-m_n}
{\rm
Let ${\mathcal E}_n$ denote the $\R$-algebra of $C^\infty$-functions-germs $(\R^n, 0) \to \R$. 
Then ${\mathcal E}_n$ is a local $\R$-algebra with the unique maximal ideal 
${\mathfrak m}_n = \{ h \in {\mathcal E}_n \mid h(0) = 0\}$. 
}
\ene

\bef
{\rm 
(\cite{Malgrange}) 
A local $\R$-algebra $A$ is called a {\bf differentiable algebra} if a surjective $\R$-algebra homomorphism, 
mapping $1$ to $1$, $\pi : {\mathcal E}_m \to A$, for some $m \in \NN$ is endowed. 
}
\enf

A differentiable algebra $A$ has the unique maximal ideal ${\mathfrak m}_A = \pi({\mathfrak m}_m)$.

\

Let $A$ and $B$ be differentiable algebras with the surjective homomorphisms $\pi : {\mathcal E}_m \to A$ and 
$\psi : {\mathcal E}_n \to B$ respectively. An $\R$-algebra homomorphism 
$u : A \to B$ is called a {\bf morphism} of differentiable algebras if 
there exists a $C^\infty$ map-germ $g : (\R^n, 0) \to (\R^m, 0)$ such that the diagram
$$
\begin{array}{ccc}
{\mathcal E}_m & \stackrel{g^*}{\rightarrow} & {\mathcal E}_n
\\
\pi\downarrow & & \downarrow \psi
\\
A & \stackrel{u}{\rightarrow} & B
\end{array}
$$
commutes.

\

A morphism $u : A \to B$ of differentiable algebras is called {\bf finite} (resp. {\bf quasi-finite}) if 
$B$ is a finite $A$-module via $u$ (resp. $B/{\mathfrak m}_AB$ is a finite dimensional $\R$-vector space). 

If $u$ is finite, then it is quasi-finite. Then we have: 

\bet
\label{M-preparation-theorem-for-differentiable-algebras}
{\rm ({\bf Malgrange preparation theorem on differentiable algebras}. Theorem 4.1 in \cite{Malgrange} p.73) }
Let $u : A \to B$ be a morphism of differentiable algebras. 
Then $u$ is finite if and only if it is quasi-finite. 
Moreover $b_1, \dots, b_r \in B$ generate $B$ over $A$ via $u$ if and only if 
$\overline{b}_1, \dots, \overline{b}_r \in B/{\mathfrak m}_AB$ generate $B/{\mathfrak m}_AB$ over $\R$. 
\ent

\bet
\label{M-M-preparation-theorem}
{\rm ({\bf Malgrange-Mather's preparation theorem}: Theorem 6.5, Corollary 6.6 in \cite{Brocker})}
Let $f : (\R^n, 0) \to (\R^m, 0)$ be a $C^\infty$ map-germ with the induced 
homomorphism $f^* : {\mathcal E}_m \to {\mathcal E}_n$. Let $C$ be a finite 
${\mathcal E}_n$-module. Then $C$ is a finite ${\mathcal E}_m$-module via $f^*$ if and only if 
$C/{\mathfrak m}_mC$ is a finite dimensional $\R$-vector space. 
Moreover $c_1, \dots, c_r \in C$ generate $C$ over ${\mathcal E}_m$ via $f^*$ 
if and only if $\overline{c}_1, \dots, \overline{c}_r \in C/{\mathfrak m}_mC$ generate $C/{\mathfrak m}_mC$ over $\R$. 
\ent

\noindent
{\it Proof that Theorem \ref{M-M-preparation-theorem} implies Theorem \ref{M-preparation-theorem-for-differentiable-algebras}:}
\\
Let $u$ is quasi-finite. Suppose $b_1, \dots, b_r \in B$ and 
$\overline{b}_1, \dots, \overline{b}_r \in B/{\mathfrak m}_AB$ generate $B/{\mathfrak m}_AB$ over $\R$. 
Let $g : (\R^n, 0) \to (\R^m, 0)$ and $g^* : {\mathcal E}_m \to {\mathcal E}_n$ cover 
$u : A \to B$. 
Note $B$ is a finite ${\mathcal E}_n$-module via $\psi$. In fact $1 \in B$ generates $B$ over ${\mathcal E}_n$ via the 
surjection $\psi$. Also 
note that ${\mathfrak m}_mB = \pi({\mathfrak m}_m)B \subseteq {\mathfrak m}_AB$. 
Then $\overline{b}_1, \dots, \overline{b}_r \in B/{\mathfrak m}_mB$ generate $B/{\mathfrak m}_mB$ over $\R$ via 
$u\circ \pi = \psi\circ g^*$. 
Therefore, by Theorem \ref{M-M-preparation-theorem}, $b_1, \dots, b_r \in B$ generate $B$ over $A$. 
Thus $u$ is finite. This implies also remaining statement naturally. 
\QED

\

\bef
{\rm 
A commutative ring $A$ is called a {\bf $C^\infty$-ring} if the following conditions are satisfied: 
\\
(1) $A$ contains the field $\R$ of real numbers. 
\\
(2) For any positive integer $r$, for any $a_1, \dots, a_r \in A$, and for any 
$C^\infty$ function $f \in C^\infty(\R^r)$, an element $f(a_1, \dots, a_r) \in A$ is assigned, such that 
the equality 
$$
(g(f_1, \dots, f_s))(a_1, \dots, a_r) = g(f_1(a_1, \dots, a_r), \dots, f_s(a_1, \dots, a_r))
$$
holds for any $g \in C^\infty(\R^s), f_1, \dots, f_s \in C^\infty(\R^n)$. 
\\
(3) The operations on $A$ by $C^\infty$ functions are compatible with the structure of $\R$-algebra on 
$A$, i.e. if $f$ is a polynomial, $f = P(x_1, \dots, x_r) \in \R[x_1, \dots, x_r] \subset C^\infty(\R^r)$, then 
$f(a_1, \dots, a_r)$ is equal to the element $P(a_1, \dots, a_r)$ obtained just by substitutions (see \cite{Ishikawa83}). 
}
\enf

Note that by the condition (1), a $C^\infty$-ring is naturally an $\R$-algebra. 
A $C^\infty$-ring $A$ is called a local $C^\infty$-ring if $A$ is a local $\R$-algebra. 
Let ${\mathfrak m}_A$ denote the unique maximal ideal of a local $C^\infty$-ring. 
Let $A$ be a $C^\infty$-ring. We say that $a_1, \dots, a_n \in A$ generate $A$ as the $C^\infty$-ring if 
for any $a \in A$, there exists $f \in C^\infty(\R^n)$ such that $a = f(a_1, \dots, a_n)$. $A$ is called a {\bf finitely generated} 
$C^\infty$-ring 
if there exists a finite number of elements generating $A$ as the $C^\infty$-ring. 
Let $\pi : A \to A/{\mathfrak m}_A$ denote the natural projection and $i :\R \to A$ the inclusion.

\bel
\label{differentiable-algebra}
Let $A$ be a differentiable algebra with a surjective $\R$-algebra homomorphism $\pi : {\mathcal E}_m \to A$. 
Then we have: 
\\
{\rm (1)} $A$ has the induced structure of a local $C^\infty$-ring.
\\
{\rm (2)} $A$ is generated by $\pi(x_1), \dots, \pi(x_m)$ as the $C^\infty$-ring. 
Here $(x_1, \dots, x_m)$ is a system of coordinates of $(\R^m, 0)$ centred at $0$. 
\\
{\rm (3)} $\pi\circ i : \R \to A/{\mathfrak m}_A$ is a bijection. 
\enl

\Proof
(1) 
For any positive integer $r$, for any $a_1, \dots, a_r \in A$, and for any 
$C^\infty$ function $f \in C^\infty(\R^r)$, we take a system of 
lifts $\widetilde{a}_1, \dots, \widetilde{a}_r \in {\mathcal E}_m$ for $\pi$ and 
define $f(a_1, \dots, a_r) := \pi(f(\widetilde{a}_1, \dots, \widetilde{a}_r))$. 
If we take another system of lifts 
$\widehat{a}_1, \dots, \widehat{a}_r \in {\mathcal E}_m$ for $\pi$, we have 
$$
f(\widetilde{a}_1, \dots, \widetilde{a}_r) - f(\widehat{a}_1, \dots, \widehat{a}_r) 
= \sum_{i=1}^r g_i(\widetilde{a}, \widehat{a})(\widetilde{a}_i - \widehat{a}_i) \in \Ker(\pi), 
$$
for some $C^\infty$ functions $g_i(\widetilde{x}_1, \dots, \widetilde{x}_r; \widehat{x}_1, \dots, \widehat{x}_r) \in C^\infty(\R^{2r}), 1 \leq i \leq r$. 
Thus $\pi(f(\widehat{a}_1, \dots, \widehat{a}_r)) = \pi(f(\widetilde{a}_1, \dots, \widetilde{a}_r))$. 
Moreover, take any $a \in A$. Then there exists $h \in {\mathcal E}_m$ such that $a = \pi(h)$. 
(2) 
Take an $H \in C^\infty(\R^m)$ having $h$ as the germ at $0$. Then 
$H(\pi(x_1), \dots, \pi(x_m)) = \pi(H(x_1, \dots, x_m)) = \pi(h) = a$. 
(3) 
$\R \cong {\mathcal E}_m/{\mathfrak m}_m \cong A/{\mathfrak m}_A$. 
\QED

\bep
Let $(A, {\mathfrak m}_A)$ be a local $C^\infty$-ring. 
Then the following conditions are equivalent:
\\
{\rm (1)} $A$ is finitely generated as $C^\infty$-ring and the natural map $\pi\circ i : \R \to A/{\mathfrak m}_A$ is bijective. 
\\
{\rm (2)} $A$ is a differentiable algebra in the sense of Malgrange. 
\enp

\Proof
(1) $\Rightarrow$ (2): 
Let $a_1, \dots, a_m$ be a system of generators of $A$ as $C^\infty$-ring. 
Define $\Pi : C^\infty(\R^m) \to A$ by $\Pi(f) = f(a_1, \dots, a_m)$. Then $\Pi$ is surjective. 
Set $I = \Pi^{-1}({\mathfrak m}_A)$ which is a maximal ideal of $C^\infty(\R^m)$ with 
$C^\infty(\R^m)/I \cong \R$. Then there exists a point $p \in \R^m$ such that 
$I = \{ f \in C^\infty(\R^m) \mid f(p) = 0\}$ (see Proposition 2.1 \cite{AS} for instance). 
Moreover $\Ker(\Pi) \subset I$. 
Set $J = \{ f \in C^\infty(\R^m) \mid {\mbox{\rm  the germ of \ }} f {\mbox{\rm \ at\ }} p {\mbox{\rm \ is  zero}}\}$. 
We show that $J \subseteq \Ker(\Pi)$. 
Let $h \in J$. Then there exists $k \in C^\infty(\R^m)$ such that $k(p) \not= 0$ and $hk = 0$. 
Then $0 = (hk)(a_1, \dots, a_m) = h(a_1, \dots, a_m)k(a_1, \dots, a_m)$. On the other hand 
$k(a_1, \dots, a_m) \not\in {\mathfrak m}_A$. Hence $k(a_1, \dots, a_m)$ is invertible. 
Then $\Pi(h) = (h(a_1, \dots, a_m) = 0$ and thus $h \in \Ker(\Pi)$. 
Now $\Pi : C^\infty(\R^n) \to A$ induces a surjective homomorphism $\pi' : {\mathcal E}_{\R^m, p} \to A$. 
Define $\pi : {\mathcal E}_m = {\mathcal E}_{\R^m, 0} \to A$ by $\pi(f) = \pi'(\widetilde{f})$, 
where $\widetilde{f}(x) = f(x - p)$. 
\\
The implication (2) $\Rightarrow$ (1) follows by Lemma \ref{differentiable-algebra}. 
\QED

\bee
{\rm 
Let $I = \{ h \in C^\infty(\R) \mid \exists n_0 \in \NN, h(n) = 0 (n \in \NN, n \geq n_0) \}$. Then $I$ is a maximal ideal of $C^\infty(\R)$. 
Let $A = C^\infty(\R)_I$ be the localisation (a localisation at infinity). 
Then $A$ is an $\R$-algebra with the unique maximal ideal 
${\mathfrak m}_A$. However $A$ is not a differentiable algebra in the sense of Malgrange. 
In fact, the quotient field $A/{\mathfrak m}_A \cong C^\infty(\R)/I$ is a Robinson's {\it hyper-real number field} 
\cite{Robinson}. 
}
\ene

We call an $\R$-algebra homomorphism $u$ a {\bf $C^\infty$-ring homomorphism} if 
$$
u(f(a_1, \dots, a_r)) = f(u(a_1), \dots, u(a_r)), 
$$
for any $r \geq 1$, for any $a_1, \dots, a_r \in A$ and for any $f \in C^\infty(\R^r)$. 

\bel
Let $\varphi : {\mathcal E}_m \to {\mathcal E}_n$ be an $\R$-algebra homomorphism.
Then the following conditions are equivalent: 
\\
{\rm (1)} There exists a $C^\infty$ map-germ $g : (\R^n, 0) \to (\R^m, 0)$ such that $\varphi = g^*$. 
\\
{\rm (2)} $\varphi$ is a $C^\infty$-ring homomorphism. 
\enl

\Proof
(1) $\Rightarrow$ (2): Let $a_1, \dots, a_r \in {\mathcal E}_m$ and $h \in C^\infty(\R^r)$. 
Then 
$$
h(\varphi(a_1), \dots, \varphi(a_r)) = h(g^*a_1, \dots, g^*a_r) = h\circ (a_1, \dots, a_r)\circ g 
= g^*(h(a_1, \dots, a_r)) = \varphi(h(a_1, \dots, a_r)). 
$$
(2) $\Rightarrow$ (1): Let $x_1, \dots, x_m$ be coordinates of $(\R^m, 0)$. Then $\varphi(x_1), \dots, \varphi(x_m) \in {\mathfrak m}_n$. 
Take representatives $\widetilde{g}_i : U \to \R$ of $\varphi(y_i)$ over a common open neighbourhood of $0$ in $\R^n$, $(1 \leq i \leq m)$. 
We set $\widetilde{g} = (\widetilde{g}_1, \dots, \widetilde{g}_m) : U \to \R^m$. Then $\widetilde{g}(0) = 0$. Take the germ $g : (\R^n, 0) 
\to (\R^m, 0)$ of $\widetilde{g}$ at $0$. Let $h \in {\mathcal E}_m$. Take a representative ${\widetilde h} \in C^\infty(\R^m)$. 
Then we have 
$\varphi(h) = \varphi(\widetilde{h}(x_1, \dots, x_m)) = \widetilde{h}(\varphi(x_1), \dots, \varphi(x_m)) = h\circ g = g^*(h)$. 
Therefore $\varphi = g^*$. 
\QED

\bel
Let $u : A \to B$ be an $\R$-algebra homomorphism of differentiable algebras. 
Then the following conditions are equivalent:
\\
{\rm (1)} $u$ is a morphism of differentiable algebras. 
\\
{\rm (2)} $u$ is a $C^\infty$-ring homomorphism. 
\enl

\Proof
(1) $\Rightarrow$ (2): 
Let $a_1, \dots, a_r \in A$ and $f \in C^\infty(\R^r)$. Take $\widetilde{a}_i \in {\mathcal E}_m$ with $\pi(\widetilde{a}_i) = a_i$. 
Then $\psi(g^*\widetilde{a}_i) = u(a_i)$. 
Then 
$u(f(a_1, \dots, a_r) = u(\pi(f(\widetilde{a}_1, \dots, \widetilde{a}_r))) = \psi(f\circ(\widetilde{a}_1, \dots, \widetilde{a}_r)\circ g) 
= \psi(f(g^*\widetilde{a}_1, \dots, g^*\widetilde{a}_r)) = f(u(a_1), \dots, u(a_r))$. 
\\
(2) $\Rightarrow$ (1): 
Take $g_i \in {\mathcal E}_n$ with $u(\pi(x_i)) = \psi(g_i)$. Since $\psi(g_i) \in {\mathfrak m}_B$, we have $g_i \in {\mathfrak m}_m$. 
Set $g = (g_1, \dots, g_m) : (\R^n, 0) \to (\R^m, 0)$. 
Let $h \in {\mathcal E}_m$ and take a representative 
$H \in C^\infty(\R^m)$ of the germ $h$. Then we have 
$u(\pi(h)) = u(H(\pi(x_1), \dots, \pi(x_m))) = H(u(\pi(x_1)), \dots, u(\pi(x_m))) = H(\psi(g_1), \dots, \psi(g_m)) = \psi(H(g_1, \dots, g_m) 
= \psi(g^*(h))$. 
\QED

\

\begin{flushleft}
Goo ISHIKAWA, \\
Department of Mathematics, Hokkaido University, 
Sapporo 060-0810, Japan. \\
e-mail : ishikawa@math.sci.hokudai.ac.jp \\
\end{flushleft}


{\footnotesize

\begin{thebibliography}{999}

%
%
%

\bibitem{Agrachev}
A.A. Agrachev, {\it Rolling balls and octonions}, 
Proc. Steklov Inst. of Math, {\bf 258} (2007), 13--22. 

\bibitem{AS}
A.A. Agrachev, Y.L. Sachkov, {\it Control Theory from the Geometric Viewpoint, } 
Springer-Verlag, Berlin Heidelberg (2004). 

\bibitem{AG}
M.A. Akivis, V.V. Goldberg, 
{\it 
Differential geometry of varieties with degenerate Gauss maps, } 
CMS Books in Mathematics, 
{\bf 18}, Springer-Verlag, New York, (2004).

\bibitem{AN}
D. An, P. Nurowski, {\it Twistor space for rolling bodies,}
Commun. Math. Phys. {\bf 326} (2014), 393--414. 


\bibitem{Arnold76}
V.I. Arnold, 
{\it 
Wave front evolution and equivariant Morse lemma, }
Comm. Pure and Appl. Math., {\bf 29} (1976) 557-582.





\bibitem{Arnold81}
V.I. Arnol'd, 
{\it Lagrangian manifold singularities, asymptotic rays and the open swallowtail, }
Funct. Anal. Appl., {\bf 15} (1981). 235--246. 


\bibitem{Arnold90}
V.I. Arnold,
{\it
Singularities of Caustics and Wave Fronts, }
Mathematics and its applications (Soviet series), {\bf 62},
Kluwer Academic Publishers., Dordrecht, (1990).


%
%
%
%



%
%
%


\bibitem{BH}
J. C. Baez, J. Huerta, 
{\it $G_2$ and the rolling ball, }
Trans. Amer. Math. Soc. {\bf 366} (2014), 5257--5293. 

\bibitem{BM}
G. Bor, R. Montgomery, 
{\it 
$G_2$ and the rolling distributions, }
Enseign. Math. {\bf 55} (2009), 157--196. 


\bibitem{Brander1}
D. Brander, 
{\it 
Pseudospherical frontals and their singularities, 
}
arXiv:1502.04876 [math.DG]. 

\bibitem{Brander2}
D. Brander, 
{\it 
Pseudospherical surfaces with singularities, 
}
arXiv:1502.04876v3 [math.DG]. 




\bibitem{Bredon}
G.E. Bredon, {\it Sheaf theory,} 
Graduate Texts in Math., {\bf 170} (2nd ed.), Springer-Verlag (1997). 

\bibitem{Bruce}
J.W. Bruce, 
{\it 
Envelopes, duality and contact structures, 
}
Proc. Symp. Pure Math., {\bf 40--1} (1983), 195--202. 


\bibitem{BG}
J.W. Bruce, P.J. Giblin, 
{\it 
Curves and Singularities, 
}
Cambridge Univ. Press, (1984). 


%
\bibitem{Brocker}
Th. Br\"{o}cker, 
{\it 
Differentiable Germs and Catastrophes, 
}
London Math. Soc. Lecture Note Series {\bf 17}, Cambridge Univ. Press (1975). 
%


%
%

%
%
%
%


%
%
%
%
%
%

\bibitem{CM}
A.L. Castro, R. Montgomery, 
{\it Spatial curve singularities and the monster/semple tower, }
Israel Journal of Math., {\bf 192} (2012), 381--427. 

\bibitem{Cayley}
A. Cayley, M{\' e}moire sur les coubes {\` a} double courbure et les surfaces d{\' e}veloppables, Journal de Mathematique Pure et Appliquees (Liouville), {\bf 10} (1845), 245--250 = The Collected Mathematical Papers vol. I, pp. 207--211. 

\bibitem{Chen-Izumiya}
L. Chen, S. Izumiya, 
{\it A mandala of Legendrian dualities for pseudo-spheres in 
semi-Euclidean space, 
}
Proc. Japan Acad., {\bf 85} Ser.A, (2009), 49--54. 


%
%

\bibitem{Chino-Izumiya}
S. Chino, S. Izumiya, 
{\it 
Lightlike developables in Minkowski 3-space, 
}
Demonstratio Mathematica {\bf 43--2} (2010), 387--399. 

\bibitem{Cleave}
J.P. Cleave, 
{\it The form of the tangent-developable at points of zero torsion on space curves, 
}
Math. Proc. Cambridge Philos. Soc. {\bf 88--3} (1980), 403--407.

%
%
%
%
%
%
%
%

\bibitem{FT1}
T. Fukunaga, M. Takahashi, 
{\it 
Existence and uniqueness for Legendre curves, }
Journal of Geometry, {\bf 104} (2013), 297--307. 

\bibitem{FT2}
T. Fukunaga, M. Takahashi, 
{\it 
Evolutes of fronts in the Euclidean plane, }
Journal of Singularities, {\bf 10} (2014), 92--107. 

\bibitem{FT3}
T. Fukunaga, M. Takahashi, 
{\it 
Involutes of fronts in the Euclidean plane,}
Beitr{\" a}ge zur Algebra und Geometrie, to appear, 
DOI: 10.1007/s13366-015-0275-1.

\bibitem{FT4}
T. Fukunaga, M. Takahashi, 
{\it 
Evolutes and involutes of frontals in the Euclidean plane,}
Demonstratio Mathematica, {\bf 48} (2015), 147--166. 



\bibitem{FSUY}
S. Fujimori, K. Saji, M. Umehara, K. Yamada, 
{\it Singularities of maximal surfaces, }
Math. Z. {\bf 259} (2008), 827--848. 

%
%

%
%
%
\bibitem{Givental86}
A.B. Givental, 
{\it Lagrangian imbeddings of surfaces and the open Whitney umbrella. }
Funk. Anal. i Prilozhen, {\bf 20--3} (1986), 35--41. 

%
%
%
%
%
%
%

%

\bibitem{Harvey}
F. R. Harvey, 
{\it 
Spinors and Calibration, 
}
Academic Press (1990). 
%
%





%
%
%

\bibitem{HKS}
A. Honda, M. Koiso, K. Saji, 
{\it 
Fold singularities on spacelike CMC surfaces in Lorentz-Minkowski space, } 
arXiv:1509.03050 [math.DG], to appear in Hokkaido Mathematical Journal. 


\bibitem{Ishikawa83}
G. Ishikawa,
{\it
Families of functions dominated by distributions of $C$-classes
of mappings,
}
Ann. Inst. Fourier {\bf 33--2} (1983), 199--217.


\bibitem{Ishikawa92}
G. Ishikawa, 
{\it 
Parametrization of a singular Lagrangian variety, 
}
Trans. Amer. Math. Soc., {\bf 331--2} (1992), 787--798. 


\bibitem{Ishikawa92-2}
G.Ishikawa, 
{\it 
The local model of an isotropic map-germ arising from one dimensional symplectic 
reduction,} 
Math. Proc. Camb. Philo. Soc., {\bf 111--1} (1992), 103--112. 


\bibitem{Ishikawa93}
G. Ishikawa, 
{\it 
Determinacy of the envelope of the osculating hyperplanes to a curve, 
}
Bull. London Math. Soc., {\bf 25}(1993), 603--610. 


\bibitem{Ishikawa94}
G. Ishikawa, 
{\it Parametrized Legendre and Lagrange varieties, 
}
Kodai Math. J., 
{\bf 17--3} (1994, October), pp.442--451. 


\bibitem{Ishikawa95}
G. Ishikawa, 
{\it 
Developable of a curve and determinacy relative to osculation-type, }
Quart. J. Math. Oxford, {\bf 46} (1995) 437--451. 

\bibitem{Ishikawa96} 
G. Ishikawa,
{\it Symplectic and Lagrange stabilities of open Whitney umbrellas, }
Invent. math., {\bf 126-2} (1996), 215--234.



\bibitem{Ishikawa99}
G. Ishikawa, 
{\it 
Singularities of developable surfaces, 
}
London Math. Soc. Lect. Notes Series, {\bf 263} (1999), 403--418. 


\bibitem{Ishikawa00}
G. Ishikawa, 
{\it 
Topological classification of the tangent developables of space curves, 
}
J. London Math. Soc., {\bf 62-2} (2000), 583--598. 

\bibitem{Ishikawa00-2}
G. Ishikawa, 
{\it Several questions on singularities: theories and applications,} 
RIMS K{\= o}ky{\= u}roku, {\bf 1122} (2000), 35--48. 

\bibitem{Ishikawa04} 
G. Ishikawa, 
{\it Classifying singular Legendre curves by contactomorphisms, }
J. of Geom. Physics, {\bf 52--2} (2004), 113-126. 

%



\bibitem{Ishikawa05}
G. Ishikawa, 
{\it 
Infinitesimal deformations and stability of singular Legendre submanifolds, 
}
Asian J. Math., {\bf 9--1} (2005), 133--166. 


%

%


\bibitem{Ishikawa10}
G. Ishikawa, 
{\it 
Singularities of flat extensions 
from generic surfaces with boundaries,  
}
Differential Geometry and its Applications {\bf 28} (2010), 341--354. 

\bibitem{Ishikawa12}
G. Ishikawa, 
{\it 
Generic bifurcations of framed curves in a space form and their envelopes, 
}
Topology and its Appl., {\bf 159} (2012), 492--500. 


\bibitem{Ishikawa12-2}
G. Ishikawa, 
{\it Singularities of tangent varieties to curves and surfaces,} 
Journal of Singularities, {\bf 6} (2012), 54--83. 

\bibitem{Ishikawa13}
G. Ishikawa, 
{\it Tangent varieties and openings of map-germs, } 
RIMS K\={o}ky\={u}roku Bessatsu, {\bf B38} (2013), 119--137. 

\bibitem{Ishikawa14}
G. Ishikawa, 
{\it Openings of differentiable map-germs and unfoldings, } 
Topics on Real and Complex Singularities, 
Proceedings of the 4th Japanese-Australian Workshop (JARCS4), Kobe 2011, World Scientific (2014), 87--113. 

\bibitem{Ishikawa15}
G. Ishikawa, 
{\it Classification problems on singularities of mappings and their applications, }
Sugaku Exposition, {\bf 28-2} (2015), 189--214.

\bibitem{Ishikawa15-2}
G. Ishikawa, 
{\it Singularities of Curves and Surfaces in Various Geometric Problems, } 
CAS Lecture Notes 10, Exact Sciences, Warsaw University of Technology (2015). 

\bibitem{IJ03}
G. Ishikawa, S. Janeczko, 
{\it 
Symplectic bifurcations of plane curves and 
isotropic liftings, }
Quart. J. Math., {\bf 54} (2003), 1--30. 

%
%

\bibitem{IJ08}
G. Ishikawa, S. Janeczko,  
{\it 
Bifurcations in symplectic space, } 
Banach Center Publ., {\bf 82} (2008), 111-124. 
%
%
%

\bibitem{IKY1}
G. Ishikawa, Y. Kitagawa, W. Yukuno, 
{\it Duality of singular paths for (2,3,5)-distributions, } 
Journal of Dynamical and Control Systems, {\bf 21} (2015), 155--171. 

\bibitem{IKY2}
G. Ishikawa, Y. Kitagawa, W. Yukuno, 
{\it Duality on geodesics of Cartan distributions and sub-Riemannian pseudo-product structures, }Demonstratio Mathematica. {\bf 48-2} (2015), 193--216． 

\bibitem{IM}
G. Ishikawa, Y. Machida, 
{\it 
Singularities of improper affine spheres and 
surfaces of constant Gaussian curvature, 
}
Intern. J. of Math., {\bf 17--3} (2006), 269--293. 

\bibitem{IM2}
G. Ishikawa, Y. Machida, 
{\it 
Monge-Amp{\`e}re systems with Lagrangian pairs, 
}
Symmetry, Integrability and Geometry: 
Methods and Applications (SIGMA), 11 (2015), 081, 32 pages. 


\bibitem{IMT1}
G. Ishikawa, Y. Machida, M. Takahashi, 
{\it 
Asymmetry in singularities of tangent surfaces in contact-cone Legendre-null duality, 
}
Journal of Singularities, {\bf 3} (2011), 126--143. 

\bibitem{IMT2}
G. Ishikawa, Y. Machida, M. Takahashi, 
{\it 
Singularities of tangent surfaces  
in Cartan's split $G_2$-geometry, 
}
Asian J. of Math., {\bf 20--2} (2016), 353--382. 


\bibitem{IMT3}
G. Ishikawa, Y. Machida, M. Takahashi, 
{\it 
Geometry of $D_4$ conformal triality
and singularities of tangent surfaces, 
}
J. of Singularities, {\bf 12} (2015), 27--52. 


\bibitem{IMT4}
G. Ishikawa, Y. Machida, M. Takahashi, 
{\it 
$D_n$-geometry and singularities of tangent surfaces, 
}
Hokkaido University Preprint Series in Mathematics {\#}1058, (2014). 
To appear in RIMS K\={o}ky\={u}roku Bessatsu (2016). 

\bibitem{IMo}
G. Ishikawa, T. Morimoto, 
{\it 
Solution surfaces of Monge-Amp\`{e}re equations, }
Diff. Geom. Appl., {\bf 14--2} (2001), 113--124.

\bibitem{IY}
G. Ishikawa, T. Yamashita, 
{\it Affine connections and singularities 
of tangent surfaces to space curves, }
arXiv:1501.07341 [math.DG].

\bibitem{IY1}
G. Ishikawa, T. Yamashita, 
{\it 
Singularities of tangent surfaces to generic space curves, }
arXiv:1602.02458 [math.DG].  

\bibitem{IY2}
G. Ishikawa, T. Yamashita, 
{\it Singularities of tangent surfaces to directed curves, }
in preparation. 

%

\bibitem{Izumiya}
S. Izumiya, 
{\it 
Legendrian dualities and spacelike hypersurfaces in the lightcone, 
}
Moscow Mathematical Journal {\bf 9} (2009), 325--357. 



%
%
%

\bibitem{IFRT}
S. Izumiya, Maria del Carmen Romero Fuster, Maria Aparecida Soares Ruas, 
F. Tari, 
{\it Differential Geometry from a Singularity Theory Viewpoint, }
World Scientific Publishing Co. (2015). 


\bibitem{IKY}
S. Izumiya, H. Katsumi, T. Yamasaki, 
{\it 
The rectifying developable and the spherical Darboux image of a space curve, }
Caustics '98, Banach Center Publ., {\bf 50} (1999), 137--149.



\bibitem{INS}
S. Izumiya, T. Nagai, K. Saji, 
{\it Great circular surfaces in the three-sphere, 
}
Diff. Geom. and its Appl., {\bf 29--3} (2011), 409--425. 


\bibitem{IPS}
S. Izumiya, D. Pei, T. Sano, 
{\it Singularities of hyperbolic Gauss maps, }
Proc. London Math. Soc., 
{\bf 86} (2003), 485--512. 

\bibitem{IPS2}
S. Izumiya, D. Pei, T. Sano, 
{\it Horospherical surfaces of curves in hyperbolic space, }
Publ. Math. Debrecen, {\bf 64} (2004), 1--13. 

\bibitem{IPT}
S. Izumiya, D. Pei, M. Takahashi, 
{\it Singularities of evolutes of hypersurfaces in hyperbolic space, }
Proc. Edinburgh Math. Soc., {\bf 47} (2004), 131--153. 

\bibitem{IS}
S. Izumiya, K. Saji, 
{\it 
The mandala of Legendrian dualities for pseudo-spheres in Lorenz-Minkowski space 
and \lq\lq flat" spacelike surfaces. 
}
J. of Singularities {\bf 2} (2010), 92--127. 


\bibitem{Kabata}
Y. Kabata, {\it 
Recognition of plane-to-plane map-germs, } 
Topology and its Applications, 
{\bf 202-1} (2016), 216--238. 


%
%
%

%
%
\bibitem{KRSUY}
M. Kokubu, W. Rossman, K. Saji, M. Umehara, K. Yamada, 
{\it Singularities of flat fronts in hyperbolic space, } 
Pacific J. Math., {\bf 221--2} (2005), 303--351.
%
%
%
%
\bibitem{Kossowski04}
M. Kossowski, 
{\it Realizing a singular first fundamental form as a nonimmersed surface in 
Euclidean $3$-space, }
J. Geom. {\bf 81} (2004), 101--113. 
%
%

%
%

\bibitem{Lawrence}
S. Lawrence, 
{\it 
Developable surfaces, their history and application, 
}
Nexus Network J., {\bf 13--3} (2011), 701--714. 

\bibitem{Malgrange}
B. Malgrange, 
{\it Ideals of Differentiable Functions, }
Oxford Univ. Press (1966). 

\bibitem{Mather}
J.N. Mather, 
{\it 
Stability of $C^\infty$ mappings III: 
Finitely determined map-germs,}
Publ. Math. I.H.E.S., {\bf 35} (1968), 279--308. 




%

%
%
%
%


\bibitem{Mond1}
D. Mond, 
{\it On the tangent developable of a space curve, 
}
Math. Proc. Cambridge Philos. Soc. {\bf 91--3} (1982), 351--355.

\bibitem{Mond2}
D. Mond, 
{\it 
Singularities of the tangent developable surface of a space curve, 
}
Quart. J. Math. Oxford Ser. (2) {\bf 40} (1989), 79--91. 

\bibitem{Mond3}
D. Mond, 
{\it Deformations which preserve the non-immersive locus of a map-germ, 
}
Math. Scand., {\bf 66} (1990), 21--32. 



\bibitem{Montgomery}
R. Montgomery, 
{\it 
A tour of subriemannian geometries, their geodesics and applications, }
Mathematical Surveys and Monographs, {\bf 91}, Amer. Math. Soc., Providence, RI, (2002).

\bibitem{MZ1}
R. Montgomery, M. Zhitomirskii, 
{\it
Geometric approach to Goursat flags, 
}
Ann. de L'institut Henri Poincar\'{e}, 
Analyse non Lin\'{e}aire, {\bf 18--4} (2001), 459--493. 

\bibitem{MZ2}
R. Montgomery, M. Zhitomirskii, 
{\it Points and Curves in the Monster Tower, 
}
Memoirs of Amer. Math. Soc., {\bf 956} (2010). 
%
\bibitem{Morimoto}
T. Morimoto, 
{\it 
Monge-Amp{\`e}re equations viewed from contact geometry,
} 
Banach Center Publ., {\bf 39} 
(1998), pp. 105--121. 
%
%

\bibitem{MU}
S. Murata and M. Umehara, 
{\it Flat surfaces with singularities in Euclidean 3-space,}
J. Differential Geom., 
{\bf 82--2} (2009), 279--316.

\bibitem{Naokawa}
K. Naokawa, 
{\it Singularities of the asymptotic completion of developable M{\" o}bius strips, }
Osaka J. of Math., {\bf 50-2} (2013), 425--437. 


%
%
%

\bibitem{Nuno-Ballesteros}
J.J. Nu\~{n}o Ballesteros, 
{\it Unfolding plane curves with cusps and nodes, } 
to appear in Proc. Roy. Soc. Edinburgh Sect. A.. 

\bibitem{NS}
J.J. Nu\~{n}o Ballesteros, O. Saeki, 
{\it Singular surfaces in $3$-manifolds, 
the tangent developable of a space curve and the dual of an immersed surface in $3$-space, }
Real and complex singularities (S\~{a}o Carlos, 1994), 49--64, 
Pitman Res. Notes Math. Ser., {\bf 333}, Longman, Harlow, (1995). 


%
%
%
%

\bibitem{O'Neill}
B. O'Neill, 
{\it Semi-Riemannian Geometry With Applications to Relativity,}
Pure and Applied Mathematics {\bf 103}, Academic Press (1983). 
%
%

\bibitem{Pohl}
W.F. Pohl, 
{\it 
The self-linking number of a closed space curve, 
}
Journal of Mathematics and Mechanics, {\bf 17-10} (1968), 975--985. 


\bibitem{Porteous0}
I.R. Porteous,
{\it The normal singularities of a submanifold, }
J. of Diff. Geom., 
{\bf 5} (1971), 543--564 . 
%
%
%
\bibitem{Porteous}
I.R. Porteous,
{\it Geometric Differentiation, for the Intelligence of 
Curves and Surfaces}, 
Cambridge Univ. Press, Cambridge, (1994). 


\bibitem{Porteous2}
I.R. Porteous, 
{\it 
Clifford Algebras and the Classical Groups, } 
Cambridge Studies in Adv. Math., {\bf 50}, 
Cambridge Univ. Press (1995). 
%
%
%

\bibitem{Robinson}
A. Robinson, {\it Non-standard Analysis, } 
Princeton University Press (1996). 

%
%
%
%
%
%
\bibitem{Saji10}
K. Saji, 
{\it Criteria for singularities of smooth maps from the plane into the plane and their applications, }
Hiroshima Math. J. {\bf 40} (2010), 229--239. 

\bibitem{Saji11}
K. Saji, 
{\it 
Criteria for $D_4$ singularities of wave fronts, 
}
Tohoku Math. J.,  {\bf 63--1} (2011), 137--147. 

\bibitem{SUY09-1}
K. Saji, M. Umehara, K. Yamada, 
{\it The geometry of fronts, } Annals of Math., 
{\bf 169-2} (2009), 491--529. 

\bibitem{SUY09-2}
K. Saji, M. Umehara, K. Yamada, 
{\it $A_k$ singularities of wave fronts, }
Math. Proc. Camb. Philos. Soc. {\bf 146-3} (2009), 731--746.

\bibitem{SUY12}
K. Saji, M. Umehara, K. Yamada, 
{\it Coherent tangent bundles and Gauss-Bonnet
formulas for wave fronts, }
J. Geom.Anal., {\bf 22} (2012), 383--409. 

%


\bibitem{Shcherbak1}
O.P. Shcherbak, 
{\it 
Projective dual space curves and Legendre singularities, 
}
Trudy Tbiliss Univ., 232--233 (1982), 280--336. 


%

%






\bibitem{Thom}
R. Thom, 
{\it 
Sur la th\'{e}orie des enveloppes, 
}
Journ. de Math., {\bf 41--2} (1962), 177--192. 


%
%


\bibitem{Wall81}
C.T.C. Wall, 
{\it 
Finite determinacy of smooth map-germs, }
Bull. London Math. Soc., {\bf 13} (1981), 481--539. 


%
%
%

\bibitem{Yamaguchi08}
K. Yamaguchi, 
{\it Geometry of linear differential systems towards contact geometry of second order,}
Symmetries and Overdetermined Systems of Partial Differential Equations, 
Vol. {\bf 144} of the series The IMA Volumes in Mathematics and its Applications (2008), pp 151--203. 

%
%
%
\bibitem{Zak}
F.L. Zak, 
{\it 
Tangents and Secants of Algebraic Varieties, }
Transl. of Math. Monographs {\bf 127}, Amer. Math. Soc.,  (1993). 
%
%
\bibitem{Zakalyukin1}
V. M. Zakalyukin, 
{\it 
Lagrangian and Legendrian singularities, 
}
Funct. Anal. Appl., {\bf 10} (1976), 23-31.

\bibitem{Zakalyukin2}
V. M. Zakalyukin, 
{\it 
Reconstructions of fronts and caustics depending 
on a parameter and versality of mappings, }
J. Soviet Math., {\bf 27} (1983), 2713-2735.
%
%
%

\bibitem{ZK}
V.M. Zakalyukin, A.N. Kurbatskii, 
{\it Envelope singularities of families of planes in control theory, }
Proc. Steklov Inst. Math., {\bf 262-1} (2008), 66--79. 

%
%
%



\end{thebibliography}

}
\end{document}